\documentclass[11pt]{article}
\usepackage{amssymb,amsmath,latexsym}

\begin{document}

\newtheorem{thm}{Theorem}[section]
\newtheorem{lemma}[thm]{Lemma}
\newtheorem{defn}[thm]{Definition}
\newtheorem{prop}[thm]{Proposition}
\newtheorem{corollary}[thm]{Corollary}
\newtheorem{remark}[thm]{Remark}
\newtheorem{example}[thm]{Example}

\numberwithin{equation}{section}

\def\ee{\varepsilon}
\def\qed{{\hfill $\Box$ \bigskip}}
\def\MM{{\cal M}}
\def\NN{{\cal M}}
\def\BB{{\cal B}}
\def\LL{{\cal L}}
\def\FF{{\cal F}}
\def\EE{{\cal E}}
\def\QQ{{\cal Q}}
\def\AA{{\cal A}}
\def\CC{{\cal C}}
\def\R{{\bf R}}
\def\N{{\mathbb N}}
\def\E{{\bf E}}
\def\F{{\bf F}}
\def\H{{\bf H}}
\def\P{{\bf P}}
\def\Q{{\bf Q}}
\def\S{{\bf S}}
\def\J{{\bf J}}
\def\K{{\bf K}}
\def\F{{\bf F}}
\def\A{{\bf A}}
\def\loc{{\bf loc}}
\def\eps{\varepsilon}
\def\semi{{\bf semi}}
\def\wh{\widehat}
\def\pf{\noindent{\bf Proof.} }
\def\dim{{\rm dim}}
\def\beq{\begin{equation}}
\def\eeq{\end{equation}}

\title{\Large \bf Generalized 3G theorem and
application to relativistic stable process on non-smooth open
sets}

\author{Panki Kim\\
Department of Mathematics\\
Seoul National University\\
Seoul 151-742, Republic of Korea\\
Email: pkim@snu.ac.kr \smallskip \\
Telephone number: 82-2-880-4077\\
Fax number: 82-2-887-4694 \\
 and
\smallskip \\
 Young-Ran Lee  \\
Department of Mathematics\\
University of Illinois \\
Urbana, IL 61801, USA\\
Email:yrlee4@math.uiuc.edu \smallskip \\
Telephone number: 1-217-265-6756\\
Fax number:1-217-333-9576  }

\maketitle

\newpage

\vspace*{1.5truein}

\begin{abstract}
Let $G(x,y)$ and $G_D(x,y)$ be the Green functions of  rotationally
invariant symmetric
$\alpha$-stable process
 in $\R^d$ and in an open set $D$ respectively, where $0<\alpha < 2$.
The inequality $ G_D(x,y)G_D(y,z)/G_D(x,z)\le c(G(x,y)+G(y,z))$ is a
very useful tool in studying (local) Schr\"{o}dinger operators. When
the above inequality is true with $c=c(D) \in (0,\infty)$, then we
say that the 3G theorem holds in $D$.

In this paper, we establish a generalized version of 3G theorem when
$D$ is a bounded $\kappa$-fat open set, which includes a bounded
John domain. The 3G we consider is of the form
$G_D(x,y)G_D(z,w)/G_D(x,w)$, where $y$ may be different from $z$.
When $y=z$, we recover the usual 3G.

The 3G form $G_D(x,y)G_D(z,w)/G_D(x,w)$ appears in non-local
Schr\"{o}dinger operator theory. Using our generalized 3G theorem,
we give a concrete class of functions belonging to the non-local
Kato class, introduced by Chen and Song, on $\kappa$-fat open sets.

As an application, we discuss relativistic $\alpha$-stable processes
(relativistic Hamiltonian when $\alpha=1$) in $\kappa$-fat open
sets. We identify the Martin boundary and the minimal Martin
boundary with the Euclidean boundary for relativistic
$\alpha$-stable processes in $\kappa$-fat open sets. Furthermore, we
show that relative Fatou type theorem is true for relativistic
stable processes in $\kappa$-fat open sets.

The main results of this paper hold for a large class of symmetric Markov processes,
as are  illustrated in the last section of this paper. We also discuss the generalized 3G
theorem for a large class of symmetric stable  L\'{e}vy processes.

\end{abstract}

\vspace{.6truein}

\bigskip

\vfill \eject

\begin{doublespace}

\section{Introduction}

The 3G theorem is crucial to prove the conditional gauge theorem
which says that the ratio of Green functions of the symmetric
$\alpha$-stable process with $0<\alpha \le 2$  and its local
perturbation in the Kato class, $K_{d,\alpha}$, is either
identically infinite or bounded above and below by two positive
numbers. When $\alpha =2$, we have a Brownian motion. Here, the Kato
class, $K_{d,\alpha},\ 0<\alpha \le 2,\ d \geq 2$, is the set of
Borel functions $q$ on $\R^d$ satisfying
$$\lim _{r \downarrow 0} \sup_{x \in \R^d} \int_{|x-y| \leq r} \frac{|q(y)|}{|x-y|^{d-\alpha}} dy=0.$$
For $d=\alpha=2$, $-\ln |x-y|$ replaces $|x-y|^{\alpha -d}$ (see
\cite{CFZ}).

Cranston, Fabes and Zhao \cite{CFZ} proved the 3G theorem of the
Brownian motion in a bounded Lipschitz domain $D \subset \R^d,\ d
\geq 3$, saying that there exists a positive constant $c$ depending
only on $D$ such that
\begin{equation} \label{3G:Brownian}
\frac{G_D(x,y) G_{D}(y,z)} { G_{D}(x,z)} \,\le \,c\,
 \{|x-y|^{2-d}+ |y-z|^{2-d}\},\qquad x,y,z \in D,
\end{equation}
where $G_D$ is the Green function of the Brownian motion in $D$.
(\ref{3G:Brownian}) was recently extended to bounded uniformly John
domains, $D \subset \R^d,\  d \geq 3$, by Aikawa and Lundh \cite{AL}
(see \cite{BaBu,HZ,S,Z} for $d=2$). The conditional gauge theorem
for the Brownian motion is of great importance to study the Green
function of the nonrelativistic Schr\"{o}dinger operator with a
local perturbation in $D$, $-\Delta +q(x)$, for $q$ in the Kato
class, $K_{d,2},\ d \geq 2$. Here, the free Hamiltonian, $-\Delta$,
is related to the kinetic energy and $q$ is its potential energy.
The point of the negative sign in front of the Laplacian is to have
that the spectrum of the free Hamiltonian is $[0, \infty)$.

For the symmetric $\alpha$-stable process with $\alpha \in (0,2)$,
the 3G theorem is also used to show the conditional gauge theorem.
The infinitesimal generator of the symmetric $\alpha$-stable process
is the fractional Laplacian $(-\Delta )^{\alpha /2}$ which is
non-local (for example, see \cite{CS1}). The 3G theorem for the
symmetric $\alpha$-stable process was
 first established and used by Chen and Song \cite{CS1, CS2}
 in a bounded $C^{1,1}$ domain $D \subset \R^d,\ d
\geq 2$. They prove that there exists a positive constant
$c=c(D,\alpha)$
 such that
\begin{equation} \label{3G:alpha}
\frac{G_D(x,y) G_{D}(y,z)} { G_{D}(x,z)}\, \le\, c\,
 \frac{|x-z|^{d-\alpha}}{|x-y|^{d-\alpha}|y-z|^{d-\alpha}},\qquad x,y,z \in
 D,
\end{equation}
where $G_D$ is the Green function of the symmetric $\alpha$-stable
process in $D$. Later, (\ref{3G:alpha}) was extended to a more
general domain which is called a bounded $\kappa$-fat open set (a
disconnected analogue of John domain, for the definition see
Definition \ref{fat}) by Song and Wu \cite{SW}.

However, to prove the conditional gauge theorem of the symmetric
$\alpha$-stable process with a certain class (the non-local Kato
class) of a non-local perturbation (for details, see section
\ref{sec:kato}), we need to generalize the 3G theorem. It is the
main objective of this paper to establish the inequality below.

\medskip

\begin{thm}[the Generalized 3G theorem]\label{t:3Gt}
Suppose that $D$ is a bounded  $\kappa$-fat open set and that $G_D(x,y)$
is the Green function for the symmetric $\alpha$-stable process in
$D$ with $\alpha \in (0,2)$ and $d \ge 2$. Then there exist positive
constants $c=c(D, \alpha)$ and $ \gamma <\alpha$ such that for every
$x, y, z, w \in D$
\begin{equation} \label{3G}
\frac{G_D(x,y) G_{D}(z,w)} { G_{D}(x,w)} \le c
\left(\frac{|x-w|\wedge |y-z|}{|x-y|  } \vee 1 \right)^\gamma
\left(\frac{|x-w|\wedge |y-z|}{|z-w|} \vee 1 \right)^\gamma
\frac{|x-w|^{d-\alpha}} {|x-y|^{d-\alpha} |z-w|^{d-\alpha}}
\end{equation}
where $a \wedge b := \min \{ a, b\}$, $a \vee b := \max \{ a, b\}$,
here and below.
\end{thm}

\medskip

The reason to call Theorem \ref{t:3Gt} the Generalized 3G theorem is
that it plays the same role in non-local perturbations as the 3G
theorem (\ref{3G:alpha})
 in local perturbations, and that we
recover the (classical) 3G theorem (\ref{3G:alpha}) by letting $y=z$ in (\ref{3G}). The
inequality (\ref{3G}) is not true without the first two factors in
the right-hand side of (\ref{3G}) (see Remark \ref{rmk2}).

In order to describe why one needs the generalized 3G theorem, we
introduce the measure $\P_x^w$ for the conditioned process obtained from
killed process $X^D$ (killed upon leaving  $D$)
through Doob's $h$-transform with $h(\cdot)=G_D(\cdot,w)$
(for details, see section 4 below) and use notations
$$
G(x,y,z,w):= \frac{G_D(x,y) G_{D}(z,w)} { G_{D}(x,w)} \quad
\mbox{and}\quad G(x,z,w):=G(x,z,z,w).
$$

When $\alpha =2$, the (killed) symmetric $\alpha$-stable process
$X^D$ in $D$ is the Brownian motion with infinitesimal generator
$-\Delta$ in $D$. For $q$ in the Kato class, $K_{d,2}$, define a
Schr\"{o}dinger semigroup by
\begin{equation}\label{e:I2}
Q_t f(x) :=\E_x \left[ \exp\left(\int_0^t q(X^D_s)ds \right)
f(X^D_t) \right], \quad t \ge 0, \, x \in D.
\end{equation}
A major tool to prove the comparison (conditional gauge theorem)
between Green functions, $G_D$ and $V_D$, of $-\Delta$ and $-\Delta
+q$, respectively,
 is 3G theorem and the identity
\begin{equation}\label{e:I1}
\E^w_x \left[ \int_0^{\tau_D} q(X^D_s)ds \right]
 = \int_D G(x,z,w) q(z) dz, \quad (x,w)
\in D \times D
\end{equation}
(see page 191 in \cite{CZ}). Then the conditional gauge theorem says
that the comparison between Green functions is given by the form:
$$
\frac{V_D(x,w)}{G_D(x,w)} =\E^w_x \left[ \exp\left( \int_0^{\tau_D} q(X^D_s)ds\right)\right]  , \quad (x,w) \in D \times D
$$
and the above is either identically infinite or bounded above and
below by positive constants (see \cite{CZ} for details).

On the other hand, for $\alpha \in (0,2)$, a symmetric
$\alpha$-stable process $X^D$ in an open set $D$ with
infinitesimal generator $(-\Delta)^{\alpha/2}$ has discontinuous
sample paths. So there is a large class of additive functionals of
$X^D$ which are not continuous. Additive functionals of the form
\[
A_{q+F}(t)=\int_0^t q(X^D_s)ds + \sum_{s\le t}F(X^D_{s-}, X^D_s)
\]
constitute an important class of discontinuous additive functionals
of $X^D$ (for example, see \cite{CS7}). Here  $q$ is a
Borel-measurable function on $D$ and $F$ is some bounded Borel
measurable function on $D \times D $ vanishing on the diagonal. As a
generalization of (\ref{e:I2}), such an additive functional defines
a Schr\"{o}dinger (Feynman-Kac) semigroup by
$$ T_t f(x) :=\E_x \left[ \exp\left( A_{q+F} (t)\right) f(X^D_t) \right], \quad t \ge 0, \, x \in D
$$
having infinitesimal generator of the form
$$
(-\Delta)^{\alpha/2} f(x) \,+\, q(x) f(x) \,+\, {\cal A} (d, \,
-\alpha) \int_D (e^{F(x,y)}-1) f(y) |x-y|^{-d-\alpha} dy
$$
where ${\cal A} (d, \, -\alpha)$ is a constant depending only on $d$
and $\alpha$ (see Theorems 4.6, 4.8 and Corollary 4.9 in
\cite{CS7}). The Schr\"{o}dinger operator of the above type has been
studied in \cite{C,CK2,CS7,CS4,CS6,K}. In particular, Chen and Song
\cite{C,CS7} introduced the various new Kato classes and proved
conditional gauge theorems for a large class of Markov processes.
The new non-local Kato class in \cite{C, CS7} is given in terms of
$G(x,y,z,w)$. The counterpart of (\ref{e:I1}) is
$$
\E^w_x \left[ A_{q+F} (\tau_D)\right] \,=\, \int_D G(x,z,w)  q(z)
dz \,+\, {\cal A} (d, \, -\alpha) \int_D \int_D G(x,y,z,w)
F(y,z)|y-z|^{-d-\alpha}dydz.
$$
So to give a concrete class of functions belonging to the non-local
Kato class, it is necessary to derive some inequality to handle
$G(x,y,z,w)$. In fact, if $D$ is smooth, by Green function estimates
(for example, see \cite{C2, CS3, CS1}), Chen and Song \cite{CS7}
give a sufficient condition for $F$ to be in the new Kato class.

In this paper, we establish the inequality to deal with the
generalized 3G form $G(x,y,z,w)$. Using this inequality, we show
that the sufficient condition given in \cite{CS7} works for any
bounded $\kappa$-fat open set (see Theorem \ref{t:A_2}).  An
interesting thing on the exponent, $\gamma$, is that once we find a
constant $\gamma$ less than $\alpha$, $\gamma$ will not affect the
concrete condition at all.

As an application of the concrete class of functions belonging to
the non-local Kato class, we consider the relativistic
$\alpha$-stable process in a bounded $\kappa$-fat open set. When
$\alpha=1$, the infinitesimal generator of this process is the free
relativistic Hamiltonian $\sqrt{-\Delta +m^{2}}-m$. Here the kinetic
energy of a relativistic particle is $\sqrt{-\Delta +m^{2}}-m $ ,
instead of $-\Delta$ for a nonrelativistic particle. The reason to
subtract the constant $m$ in the free Hamiltonian is to ensure that
its spectrum is $[0,\infty)$ (see \cite{CMS}). There exists a huge
literature on the properties of relativistic Hamiltonian (for
example, \cite{CMS,FL,H,Lieb,LY}). Recently relativistic
$\alpha$-stable process with the infinitesimal generator $(-\Delta
+m^{2/\alpha} )^{\alpha /2}-m$ has been studied by Chen and Song
\cite{CS4}, Kim \cite{K2},
 Kulczycki and Siudeja \cite{KuS}, Ryznar \cite{R}.

In this paper, we show that, in any bounded $\kappa$-fat open set, Green functions and Martin kernels
of relativistic stable process are comparable to the ones of
symmetric stable process. Using
recent results in \cite{CK2}, we identify the Martin boundary and
the minimal Martin boundary with the Euclidean boundary for
relativistic $\alpha$-stable process in $\kappa$-fat open set. In
\cite{K2}, the first named author studied the boundary behavior of
the ratio of two harmonic functions of Schr\"{o}dinger operators of
the above type. As a consequence of the results in \cite{K2} and the
generalized 3G inequality, we show that relative Fatou type theorem
is true for the relativistic stable process in a $\kappa$-fat open
set.

Our method works for other classes of symmetric  L\'{e}vy
processes. However to make the exposition  as transparent as
possible, we choose to present our results for (rotationally
invariant) symmetric stable processes. Extensions to other classes
of symmetric  L\'{e}vy processes are mentioned at the end of the
paper.

The rest of the paper is organized as follows. In section 2 we
collect basic definitions and properties of symmetric stable
processes. Section 3 contains the proof of the generalized 3G
theorem. In section 4, through the generalized 3G theorem, a
concrete condition is given for functions to be in the non-local
Kato class. Using the condition, the Martin boundary and the
boundary behavior of harmonic functions for relativistic stable
processes in bounded $\kappa$-fat open sets are studied. In section
5, some extensions of the main results of this paper to more general
symmetric Markov processes are given. As an example,  we show the
generalized 3G theorem is true for some class of symmetric stable
L\'{e}vy processes.

 In this paper, we will use the following convention: The
values of the constants $C_0$, $C_1$, $M$ and $\eps_1$ will remain
the same throughout this paper, while the values of the constants
$c, c_1, c_2, \cdots$ signify  constants whose values are
unimportant and which may change from location to location. In
this paper, we use ``$:=$" to denote a definition, which is  read
as ``is defined to be".    For any open set $U$, We denote by $\rho_U (x)$ the
distance of a point $x$ to the boundary of $U$, i.e.,
$\rho_U(x)=\text{dist} (x,\partial U)$.

\section{Preliminary}\label{s:prel}
In this section, We recall the definition of a rotationally
invariant symmetric stable process and collect its properties, which
we will use later.

Let $X= \{X_t\}$ denote a (rotationally invariant) symmetric
$\alpha$-stable process in $\R^d$ with $\alpha \in (0, \, 2)$ and
$d\geq 2$, that is, let $X_t$ be a L\'evy process whose transition
density $p(t, y-x)$ relative to the Lebesgue measure is given by the
following Fourier transform,
$$ \int_{\R^d} e^{i x\cdot \xi}  p(t, x) dx
        =e^{-t |\xi|^\alpha} .
$$

The Dirichlet form ${\cal E}$ with the domain $\cal F$ associated
with $X$ is given by
\begin{equation}\label{DF1}
{\cal E}(u, v):=\int_{\R^d} \hat{u}
(\xi)\bar{\hat{v}}(\xi)|\xi|^{\alpha} d\xi, \quad {\cal F}:= \{ u\in
L^2(\R^d):\int_{\R^d} |\hat{u} (\xi)|^2|\xi|^{\alpha}d\xi <\infty\},
\end{equation}
where $\hat{u} (\xi):=(2\pi)^{-d/2}\int_{\R^d}e^{i\xi\cdot y}u(y)dy$
(see Example 1.4.1 of \cite{FOT}). As usual, we define ${\cal
E}_1(u, v):={\cal E}(u, v)+(u, v)_{L^2(\R^d)}$ for $u, v\in
 \cal F$. Another expression for ${\cal E}$ is as follows:
$$
{\cal E}(u, v)=\frac12 {\cal A} (d, - \alpha)\int_{\R^d}\int_{\R^d}
\frac{(u(x)-u(y))(v(x)-v(y))} {|x-y|^{d+\alpha}}\, dxdy,
$$
where $ {\cal A}(d, -\alpha):= \alpha2^{\alpha-1}\pi^{-d/2}
\Gamma(\frac{d+\alpha}2) \Gamma(1-\frac{\alpha}2)^{-1}. $ Here
$\Gamma$ is the Gamma function defined by $\Gamma(\lambda):=
\int^{\infty}_0 t^{\lambda-1} e^{-t}dt$
 for every $\lambda > 0$. Let
\begin{equation}\label{J}
J(x, y):=\frac12{\cal A}(d, -\alpha)|x-y|^{-(d+\alpha)}.
\end{equation}
$J(x, y)dy$ is called the jumping measure of $X$. Let Cap be the
1-capacity associated with $X$ (or equivalently, with the Riesz
potential kernel ${\cal A}(d, -\alpha) |x-y|^{-d+\alpha}$). A
function $f$ is said to be quasi-continuous if for any $\eps >0$
there exists an open set $U$ such that Cap$(U) < \eps$ and
$f|_{U^c}$ is continuous. It is known that
 every function $f$ in $\cal F$ admits a quasi-continuous version.
For concepts and results related to Dirichlet forms, we refer our
readers to \cite{FOT}.

For any open set $D$, we use $\tau_D=\tau(D)$ to denote the first
exit time of $D$, i.e., $\tau_D=\tau(D)=\inf\{t>0: \, X_t\notin
D\}$. Given an open set $D\subset \R^d$, let
$X^D_t(\omega)=X_t(\omega)$ if $t< \tau_D(\omega)$ and set
$X^D_t(\omega)=\partial$ if $t\geq \tau_D(\omega)$, where $\partial$
is a coffin state added to $\R^d$. The process $X^D$, i.e., the
process $X$ killed upon leaving  $D$, is called the symmetric
$\alpha$-stable process in $D$. Throughout this paper, we use the
convention $f(\partial)=0$.

The Dirichlet form of $X^D$ is $({\cal E}, {\cal F}^D)$, where
\beq\label{F^D}
 {\cal F}^D:=\{u\in {\cal F}: u=0 \mbox{ on $D^c$
except for a set of zero capacity }\}
\eeq (cf. \cite{FOT}). Thus
for any $u, v\in {\cal F}^D$,
\begin{equation}\label{KDF}
{\cal E}(u, v)=\int_D\int_D(u(x)-u(y))(v(x)-v(y))J(x, y)dxdy
+\int_Du(x)v(x)\kappa_D(x)dx,
\end{equation}
where \beq\label{kappa_D} \kappa_D(x):={\cal A}(d,
-\alpha)\int_{D^c}|x-y|^{-(d+\alpha)}dy.
\eeq

Before we state more properties for a symmetric $\alpha$-stable
process, let's introduce the following definitions.

\medskip

\begin{defn}\label{def:har1} Let $D$ be an open subset of $\R^d$.
A locally integrable function $u$ defined on $\R^d$ taking values in
$(-\infty, \, \infty]$ and satisfying the condition $ \int_{\{x \in
\R^d ; |x| >1\}}|u(x)| |x|^{-(d+\alpha)} dx  <\infty $ is said to be

\begin{description}
\item{(1)}  harmonic for $X$ in $D$ if $$
\E_x\left[|u(X_{\tau_{B}})|\right] <\infty \quad \hbox{ and } \quad
u(x)= \E_x\left[u(X_{\tau_{B}})\right], \qquad x\in B, $$ for every
open set $B$ whose closure is a compact subset of $D$;
\item{(2}) regular  harmonic for $X$ in $D$ if it is harmonic for $X$ in $D$  and
for each $x \in D$, $$ u(x)= \E_x\left[u(X_{\tau_{D}})\right]; $$
\item{(3}) singular harmonic for $X$ in $D$ if it is harmonic for $X$ in $D$ and
it vanishes outside $D$.
\end{description}
\end{defn}

\medskip

Note that a harmonic function in an open set $D$ is continuous in
$D$ (see \cite{BB} for an analytic definition and its equivalence).
Also note that singular harmonic function  $u$ in $D$ is harmonic
with respect to $X^D$, i.e., $$ \E_x\left[|u(X^D_{\tau_{B}})|\right]
<\infty \quad \hbox{ and } \quad u(x)=
\E_x\left[u(X^D_{\tau_{B}})\right], \qquad x\in B, $$ for every open
set $B$ whose closure is a compact subset of $D$.

We have the Harnack Principle as follows.
\medskip

\begin{thm}\label{HP} (Bogdan \cite{Bo}) Let $x_{1}, x_{2}\in \R^d$,
$r>0$ such that $|x_{1}-x_{2}|< Lr$. Then there exists a constant
$J$ depending only on $d$ and $\alpha$, such that
$$ J^{-1}L^{-(d+\alpha)}u(x_{2})\leq u(x_{1}) \leq JL^{d+\alpha}u(x_{2})\, $$
for every nonnegative harmonic function $u$ in $B(x_{1}, r)\cup
B(x_{2},r)$. \end{thm} Here and below, $B(x,r)$ denotes the ball
centered at $x$ with radius $r$.

\medskip

We adopt the definition of a $\kappa$-fat open set from \cite{SW}.
\medskip
\begin{defn}\label{fat} Let $\kappa \in (0,1/2]$. We say that an open set
$D$ in $\R^d$ is $\kappa$-fat if there exists $R>0$ such that for
each $Q \in \partial D$ and $r \in (0, R)$, $D \cap B(Q,r)$ contains
a ball $B(A_r(Q),\kappa r)$ for some $A_r(Q) \in D$. The pair $(R,
\kappa)$ is called the characteristics of the $\kappa$-fat open set
$D$.
\end{defn}

\medskip

Note that every Lipschitz domain and {\it non-tangentially
accessible domain} defined by Jerison and Kenig in \cite{JK} are
$\kappa$-fat. Moreover, every {\it John domain} is  $\kappa$-fat
(see Lemma 6.3 in \cite{MV}). The boundary of a $\kappa$-fat open
set can be highly nonrectifiable and, in general, no regularity of
its boundary can be inferred.  Bounded $\kappa$-fat open set can
even be locally disconnected.

The boundary Harnack principle for $\kappa$-fat open sets was proved
in \cite{SW}.

\medskip

\begin{thm}\label{T:BH}
(Theorem 3.1 in \cite{SW}) Let $D$ be a $\kappa$-fat open set with
$(R, \kappa)$. Then there exists a constant $c=c(d, \alpha)
>1$ such that for any $Q \in \partial D$, $r \in (0,R)$ and
functions $u,v \ge 0$ in $\R^d$, regular harmonic in $D \cap
B(Q,2r)$, vanishing on $D^c \cap B(Q,2r)$, we have $$ c^{-1}
\kappa^{d+ \alpha} \frac{u(A_r(Q))}{v(A_r(Q))} \le \frac{u(x)}{v(x)}
\le c \kappa^{-d- \alpha}\frac{u(A_r(Q))}{v(A_r(Q))}, ~~~~~ x\in
D\cap B\left(Q,\frac{r}{2}\right). $$ \end{thm}

\medskip

It is well known that there is a positive continuous symmetric
function $G_D(x, y)$ on $(D\times D)\setminus \{x=y\}$ such that for
any Borel measurable function $f\geq 0$, $$ \E_x \left[
\int_0^{\tau_D} f(X_s) ds \right] =\int_D G_D (x, y) f(y) \, dy.
$$ We set $G_D$ equal to zero on the diagonal of $D \times D$ and
outside $D \times D$. Function $G_D(x, y)$ is called the Green
function of $X^D$, or the Green function of $X$ in $D$. For any $x
\in D$, $G_D( \, \cdot \, ,x)$ is singular harmonic in $D\setminus
\{x\}$ and regular harmonic in $D\setminus B(x, \eps)$ for every
$\eps >0$. When $D=\R^d$, it is well known that $G_{\R^d}(x,y) =
c|x-y|^{\alpha-d}$ where $c=c(d,\alpha)$ is a positive constant
depending only on $d$ and $\alpha$.

\section{Generalized 3G theorem}

In this section, we prove the objective of this paper, the
generalized 3G theorem \ref{t:3Gt}, in a bounded $\kappa$-fat open
set $D$. Let $D$ be a $\kappa$-fat open set with its characteristics
$(R, \kappa)$, $0< \kappa \le 1/2$. We fix $R$ and $\kappa$
throughout this section. We collect and prove lemmas to prove the
main theorem.

For any ball $B$ in $\R^d$ let  $P_B(x,z)$ be the Poisson kernel for
$X$ in $B$, i.e.,
$$\P_x(X_{\tau_B} \in A)=
\int_{A} P_B(x,z) dz, \quad A \subset \R^d \setminus B.
$$
It is well-known that
\begin{equation}\label{P_f}
P_{B(x_0,r)}(x,z)\,=\,c_1\,
\frac{(r^2-|x-x_0|^2)^{\frac{\alpha}2}}{(|z-x_0|^2-r^2)^{\frac{\alpha}2}}
\frac1{|x-z|^d}
\end{equation}
for some constant $c_1=c_1(d, \alpha) > 0$.

The next lemma can be proved easily using  (\ref{P_f}) (see
\cite{Bo,KS,SW}). Since knowing the value of $\gamma$ might be
interesting for readers, we put a proof here. Moreover, with an eye
towards the extensions in section 5, the following lemma will be
proved without using (\ref{P_f}). Instead, we will use the fact that
\begin{equation}\label{e:low}
P_{B(x,1)}(x,z) \,\ge\, c(d, \alpha) \,|x-z|^{-d -\alpha},\quad  z
\in \R^d \setminus \overline{B(x,1)}, \, x \in \R^d.
\end{equation}
Recall, from Definition \ref{fat} that, for $Q \in \partial D$ and
$r \in (0, R)$, there exists a point $A_r(Q)$ in $D \cap B(Q,r)$
such that $B(A_r(Q),\kappa r) \subset D \cap B(Q,r)$. Recall
$\rho_D(x)=\text{dist} (x,\partial D)$. Note that
$\kappa R \le \rho_D(A_r(Q))<r$.

\medskip

\begin{lemma}\label{l2.5} (Lemma 5 of \cite{Bo} and Lemma 3.6 of \cite{SW})
There exist positive constants $c=c(d, \alpha, \kappa)$ and
$\gamma=\gamma(d, \alpha, \kappa)< \alpha$ such that for all $Q\in
\partial D$ and $r\in (0, R)$, and function $u\ge 0$, harmonic in
$D\cap B(Q, r)$, we have
\begin{equation}\label{e:gamma}
u(A_s(Q))\ge c(s/r)^{\gamma}u(A_r(Q)), \qquad s\in (0, r).
\end{equation}
In fact, there exists a small constant $c_1=c_1(d, \alpha)>0$ so
that we can take $ \gamma= \alpha- c_1 (\ln (2/\kappa))^{-1}$.
\end{lemma}
\pf Without loss of generality, we may assume $ r=1$, $Q=0$ and
$u(A_r(Q)) =1$.  Let $ r_k:=(\kappa/2)^{k} $,  $A_k:= A_{r_k}(0) $
and $ B_k:=B(A_k, r_{k+1})$ for $k=0,1, \cdots. $

Note that the $B_k$'s are disjoint. So by the harmonicity of $u$,
we have
$$
u(A_k) \,\ge\, \sum_{l=0}^{k-1} \E_{A_k}\left[u(X_{\tau_{B_k}}):\,
X_{\tau_{B_k}} \in B_l \right]\\
\,=\, \sum_{l=0}^{k-1} \int_{B_l} P_{B_k}(A_k, z) u(z) dz.
$$
Since $B_k \subset B(A_k, 2 r_{k+1}) \subset D\cap B(0, r_k)$,  The Harnack
principle (Theorem 2.2) implies that
$$
\int_{B_l} P_{B_k}(A_k, z) u(z) dz \,\ge\, c_1\, u(A_l) \int_{B_l}
P_{B_k}(A_k, z) dz\,=\, c_1 \,u(A_l) \,\P_{A_k} (X_{\tau_{B_k}} \in B_l)
$$
for some constant $c_1=c_1(d, \alpha)$. By the scaling property of
$X$, we have that
\begin{eqnarray*}
&&\P_{A_k} (X_{\tau_{B_k}} \in B_l) \,=\,
\P_{r_{k+1}^{-1}A_k} (X_{\tau(r_{k+1}^{-1}B_k)} \in r_{k+1}^{-1}B_l)\\
&&\ge\, c\, (r_{k-l})^{-d}  \left ( \rho_{ r_{k+1}^{-1}B_l}
(r_{k+1}^{-1}A_k)\right)^{-d-\alpha}  \ge\, c_2 \, r_{k-l}^{
\alpha}
\end{eqnarray*}
for some constant $c_2=c_2(d, \alpha)$. In the first and second
inequalities above, we have used   (\ref{e:low}) and the fact that
$\rho_{B_l} (A_k) \le 2r_l$ respectively. Therefore,
$$
\left(r_k\right)^{-\alpha} u(A_k) \,\ge \, c_3 \sum_{l=0}^{k-1}
\left(r_l\right)^{-\alpha} u(A_l)
$$
for some constant $c_3=c_3(d, \alpha)$. Let $a_k :=
r_k^{-\alpha}u(A_k)$ so that $a_k \ge  c_3\sum_{l=0}^{k-1}  a_l$. By
induction, one can easily check that $ a_k  \ge c_4 (1+c_3/2)^{k} $
for some constant $c_4=c_4(d, \alpha)$. Thus, with $ \gamma = \alpha
- {\ln(1+\frac{c_3}2)} (\ln (2/\kappa))^{-1}, $ (\ref{e:gamma}) is
true for $s= r_k$. For the other values, use Harnack principle
(Theorem \ref{HP}). \qed

\medskip

The fact, $\gamma < \alpha$, is critical in the next section. As a
corollary of the boundary Harnack principle (Theorem \ref{T:BH}), we
immediately get the following Lemma.

\medskip

\begin{lemma}\label{l:Green_L}
There exists a constant $ c=c(D,\alpha) >1$ such that for every $Q
\in
\partial D$, $0<r < R$,
we have
\begin{equation}\label{e:CG_1}
\frac{ G_{D}(x,z_1)}{ G_{D} (y,z_1)} \, \le \,c\, \frac{ G_{D}
(x,z_2) }{ G_{D} (y,z_2)},
\end{equation}
when $x, y \in D \setminus \overline{B(Q, r)}$ and $z_1, z_2 \in D
\cap B(Q, r/4)$.
\end{lemma}

\medskip

Let $M:=2\kappa^{-1} $ and fix $z_0 \in D$ with $2R/M < \rho_D(z_0)
< R$ and let $\eps_1:= R /(12M)$. For $x,y \in D$, we define
$r(x,y): = \rho_D(x) \vee \rho_D(y)\vee |x-y|$ and
$$
\BB(x,y):=\left\{ A \in D:\, \rho_D(A) > \frac1{M}r(x,y), \,  |x-A|\vee
|y-A| < 5 r(x,y)  \right\}
$$
if $r(x,y) <\eps_1 $, and $\BB(x,y):=\{z_0 \}$ otherwise.

By the monotonicity of the Green function and the formula for the
Green function for balls (for example, see \cite{CS1}), there exists
a positive constant $C_0$ such that
\beq\label{ub}
G_D(x,y) \le C_0 |x-y|^{-d+\alpha},\qquad x,y \in D
\eeq
and that $ G_D(x,y) \ge G_{B(y, \rho _D(y))}(x,y) \ge C_0^{-1}
|x-y|^{-d+\alpha}\ \text{if }|x-y| \le \rho_D(y)/2 $. Let $C_1:=C_0
2^{d-\alpha} \rho_D(z_0)^{-d+\alpha}$ so that $G_D(\cdot, z_0)$ is
bounded above by $C_1$ on $ D \setminus B(z_0, \rho_D(z_0)/2)$. Now
we define
$$
g(x ):=  G_D(x, z_0) \wedge C_1.
$$
Note that if $\rho_D(z) \le 6 \eps_1$, then $|z-z_0| \ge \rho_D(z_0)
- 6 \eps_1 \ge \rho_D(z_0) /2$ since $6\eps_1<\rho_D(z_0)/4$, and
therefore $g(z )= G_D(z, z_0)$.

The lemma below immediately follows from the Harnack principle,
Theorem \ref{HP}.

\medskip

\begin{lemma}\label{G:g1}
There exists $c=c(D, \alpha)>0$ such that for every $x, y \in D$,
$$
c^{-1} \,g(A_1) \,\le \,g(A_2)\, \le\, c\, g(A_1) , \qquad A_1, A_2
\in \BB(x,y).
$$
\end{lemma}

\medskip

Using the Harnack principle and the boundary Harnack principle, the
following form of Green function estimates has been established by
several authors. (See Theorem 2.4 in \cite{H} and Theorem 1 in
\cite{J}.  Also see \cite{CK} for a different jump process and see
\cite{KS2} for non-symmetric diffusion.)

\medskip

\begin{thm}\label{t:Gest}
There exists $c=c(D, \alpha)>0$ such that for every $x, y \in D$
\begin{equation}\label{e:Gest}
c^{-1}\,\frac{g(x) g(y)}{g(A)^2} \,|x-y|^{-d+\alpha} \,\le\,
G_D(x,y) \,\le\, c\,\frac{g(x) g(y)}{g(A)^2}\, |x-y|^{-d+\alpha},
\quad  A \in \BB(x,y).
\end{equation}
\end{thm}

\medskip

Applying Theorem \ref{t:Gest}, we immediately get the following
inequality.

\begin{thm}\label{t:3G1}
There exists a constant $c=c(D, \alpha) >0$ such that for every  $x, y, z,w \in
D$ and $(A_{x,y}, A_{z,w}, A_{x,w}) \in \BB(x,y) \times \BB(z,w)
\times \BB(x,w)$,
\begin{equation}\label{3G_est} \frac{G_D(x,y) G_{D}(z,w)} { G_{D}(x,w)} \,\le\, c\,
\frac{g(y)g(z) g(A_{x,w})^2}{g(A_{x,y})^2 g(A_{z,w})^2}
\frac{|x-w|^{d-\alpha}}{|x-y|^{d-\alpha} |z-w|^{d-\alpha}}.
\end{equation}
\end{thm}

\medskip

Using the Harnack principle, the proof of the next lemma is
well-known (for example, see Lemma 6.7 in \cite{CZ}). We skip the
proof.
\begin{lemma}\label{G:g3} (Lemma 2.2 in \cite{H})
For every $c_1 >0$, there exists $c=c(D, \alpha, c_1)>0$ such that
for every $|x-y| \le c_1 ( \rho_D(x) \wedge \rho_D(y)) $, $ G_D(x,y)
\ge c |x-y|^{-d+\alpha}. $
\end{lemma}

\medskip

The following lemma is a direct application of Lemma 3.4 in
\cite{SW} and Lemma 4 in \cite{Bo}. However, for the extensions in
section 5, we prove here using Lemma \ref{l:Green_L}.

\medskip

\begin{lemma}\label{C:c_L} (Carleson's estimate)
There exists a constant $ c=c(D,\alpha) >1$ such that for every $Q
\in
\partial D$,  $0<r < \kappa R/4$, $y \in D  \setminus \overline{B(Q, 4r)}$
\begin{equation}\label{e:CG_3}
G_{D}(x, y) \, \le \,c\,  G_{D} ( A_r(Q), y), \quad x \in D \cap
B(Q, r).
\end{equation}
\end{lemma}

\pf
 By Lemma \ref{l:Green_L}, there exists a constant
$c=c(D,\alpha)>1$ such that for every $Q \in \partial D$,  $0<r <
\kappa R/4$, we have
$$ \frac{G_{D}(x, y)}{G_{D}(A_r(Q), y)} \, \le \,c\, \frac{G_{D}(x, z)}{G_{D} ( A_r(Q), z)}$$
for any $x \in D \cap B(Q, r)$ and $y,\ z \in D \setminus
\overline{B(Q, 4r)}$. Since $r<\kappa R/4$, i.e., $4r/\kappa <R$, we
can and do pick $z=A_{4r/\kappa}(Q)$. Then by (\ref{G:g3}) and the
inequality (\ref{ub}), there exists a positive constant
$c_1=c_1(D,\alpha)$ such that
$$\frac{G_{D}(x,A_{4r/\kappa}(Q) )}{G_{D} ( A_r(Q),
A_{4r/\kappa}(Q))}<c_1,$$ which immediately implies the inequality
(\ref{e:CG_3}). Indeed, by (\ref{ub}),
$$
G_D(x,A_{4r/\kappa}(Q))\,\le\, C_0\, |x-A_{4r/\kappa}(Q)|^{-d+\alpha} \,\le\,
c_2\, r^{-d+\alpha}$$
since $|x-A_{4r/\kappa}(Q)| \ge
\rho_D(A_{4r/\kappa}(Q))-\rho_D(x)\ge 4r-r=3r$. On the other hand,
since
$$|A_r(Q)-A_{4r/\kappa}(Q)| \,\le\, \frac{8r}{\kappa}\, \le\,
\frac{8}{\kappa ^2} (\rho_D(A_{4r/\kappa}(Q)) \wedge
\rho_D(A_r(Q))),$$ by Lemma \ref{G:g3}, we obtain
$$ G_D(A_{4r/\kappa}(Q),A_r(Q)) \,\ge\, c_3\,|A_r(Q)-A_{4r/\kappa}(Q)|^{-d+\alpha}
\, \ge\, c_3\,r^{-d+\alpha}.$$
 \qed

\medskip

For every $x,y \in D$, we denote $Q_x$ and $Q_y$ be  points on
$\partial D$ such that $\rho_D(x)=|x-Q_x|$ and
$\rho_D(y)=|y-Q_y|$ respectively. It is easy to check that  if $r(x,y) <
\eps_1$
\begin{equation}\label{e:AinB}
A_{r(x,y)}(Q_x),\, A_{r(x,y)}(Q_y) \,\in\, \BB(x,y).
\end{equation}
In fact, by definition of $ A_{r(x,y)}(Q_x)$,
$\rho_D(A_{r(x,y)}(Q_x)) \ge \kappa r(x,y) > r(x,y) /M$. Moreover,
$$
|x- A_{r(x,y)}(Q_x)|\, \le\, | x-Q_x|+|Q_x -  A_{r(x,y)}(Q_x)|
\,\le\, \rho_D(x) + r(x,y) \,\le\, 2 r(x,y)
$$
and $ |y- A_{r(x,y)}(Q_x)| \le |y-x| +|x- A_{r(x,y)}(Q_x)| \le 3
r(x,y). $ Thus we get (\ref{e:AinB}).

Recall the fact that  $g(z )=  G_D(z, z_0)$ if $\rho_D(z) < 6
\eps_1$.

\medskip

\begin{lemma}\label{G:g5_1}
There exists $c=c(D,\alpha)>0$ such that for every $x,y \in D$ with
$r(x,y) < \eps_1$,
\begin{equation}\label{e:gg5_1}
g(z)  \,\le\, c\, g(A_{r(x,y)}(Q_x)),  \quad z \in D\cap B(Q_x,
r(x,y))
\end{equation}
\end{lemma}
\pf
Note that for $z \in D\cap B(Q_x, 2 r(x,y))$, $\rho_D(z) < 2\eps_1$.
Thus $g(\cdot)=G_D(\cdot, z_0)$ in $D\cap B(Q_x, 2 r(x,y))$. By
Lemma \ref{C:c_L}, we have $g(z) =G_D(z, z_0) \le c \,G_D(
A_{r(x,y)}(Q_x) , z_0)=c\, g(A_{r(x,y)}(Q_x)) $ for $z \in D\cap
B(Q_x, r(x,y))$.
\qed

\medskip

\begin{lemma}\label{G:g5}
There exists $c=c(D,\alpha)>0$ such that for every $x,y \in D$
$$
g(x) \vee g(y) \,\le\, c\, g(A),\quad A \in \BB(x,y).
$$
\end{lemma}
\pf
If $r(x,y) \ge \eps_1$, then $\BB(x,y)=\{z_0\}$ and $g(x) \vee g(y)
\le C_1 =g(z_0) $.

Now assume $r(x,y) < \eps_1$. Since $ \rho_D (x) \le r(x,y) $, $x
\in D \cap B(Q_x, r(x,y))$. Thus, by Lemma \ref{G:g5_1},
$$
g(x)  \,\le\, c\, g(A_{r(x,y)}(Q_x)).
$$
Similarly, we obtain $g(y) \,\le c\, g(A_{r(x,y)}(Q_y)).$

The lemma follows from Lemma \ref{G:g1} and (\ref{e:AinB}). \qed

\medskip

\begin{lemma}\label{G:g4}
If $x,y,z \in D$ satisfy $ r(x,z) \le r(y,z)$, then there exists $c=c(D,\alpha)>0$ such that
$$
g(A_{x,y}) \,\le\, c\, g(A_{y,z})  \qquad \mbox{for every } (A_{x,y},
A_{y,z}) \in \BB(x,y) \times \BB(y,z).
$$
\end{lemma}
\pf
If $r(y,z) \ge \eps_1$, then $g(A_{y,z}) = C_1 \ge g(A_{x,y})$.
Thus, we assume $r(y,z) < \eps_1$. First, we note that
\begin{eqnarray*}
&&r(x,y) \,\le \,\rho_D (x) \vee \rho_D (y)   \vee (|x-z|+|z-y|)
\,\\
&&~~ \le\, \rho_D (x) \vee \rho_D (y)   \vee |x-z| + \rho_D (x) \vee
\rho_D (y)   \vee |z-y| \le\, 2(r(x,z) + r(y,z))\, \le\, 4 \,r(y,z).
\end{eqnarray*}
Here, the assumption $r(x,z) \le r(y,z)$ was used in the last
inequality.
 So  $A_{r(x,y)}(Q_y) \in D \cap B(Q_y, 4r(y,z)) $. Since
$g(\cdot)=G_D(\cdot,z_0)$ in $D \cap B(Q_y, 6r(y,z))$ and $6r(y,z) <
\kappa R/4$, by Lemma \ref{C:c_L},
$$g(w)=G_D(w,z_0) \le c G_D( A_{ 6r(y,z)} (Q_y),z_0)=cg( A_{ 6r(y,z)}
(Q_y)) $$
for $w \in D \cap B(Q_y, 6r(y,z))$. Here, we used
$|Q_y-z_0| \ge \rho_D(z_0) > 2R/M>24 r(y,z)$. In particular, with
$w=A_{r(x,y)}(Q_y)$,
$$g(A_{r(x,y)}(Q_y) ) \le c g( A_{ 6r(y,z)} (Q_y)). $$
  By Theorem \ref{HP},
$g( A_{ 6r(y,z)} (Q_y))= G_D(A_{ 6r(y,z)} (Q_y) ,z_0) \le c_1
G_D(A_{ r(y,z)} (Q_y) ,z_0) =c_1 g( A_{ r(y,z)} (Q_y))$ for some
$c_1=c_1(d,\alpha)$. Thus
$$
g(A_{r(x,y)}(Q_y)  ) \,\le\, c_2 \,g( A_{ r(y,z)} (Q_y))
$$
 for some
$c_2=c_2(D,\alpha)$. Using (\ref{e:AinB}) we apply Lemma \ref{G:g1}
both sides of the above functions and we have proved the lemma. \qed

\medskip

\begin{lemma}\label{g:tri3}
There exists $c=c(D,\alpha)>0$ such that for every  $x, y, z, w \in D$ and
$(A_{x,y}, A_{y,z}, A_{z,w}, A_{x,w}) \in \BB(x,y)\times \BB(y,z)
\times \BB(z,w) \times \BB(x,w) $,
\begin{equation}\label{e:tri3}
g(A_{x,w})^2 \,\le\, c\, \left(g(A_{x,y})^2 + g(A_{y,z})^2 +
g(A_{z,w})^2 \right).
\end{equation}
\end{lemma}
\pf Applying Lemma \ref{G:g4} to both cases $r(x,y) \le r(y,z)$
and $r(x,y) \ge r(y,z)$, we get
$$
g(A_{x,z})^2 \,\le\,  c_1 \left( g(A_{x,y})^2\,+\,
g(A_{y,z})^2\right).
$$
Obviously, using the above inequality twice, we get
$$
g(A_{x,w})^2 \,\le\, c_1 \, \left( g(A_{x,z})^2 \,+\, g(A_{z,w})^2
\right) \,\le\, c_1 \, \left( c_1 g(A_{x,y})^2 +  c_1
g(A_{y,z})^2 + g(A_{z,w})^2 \right).
$$\qed

\medskip

Recall $\gamma$ from Lemma \ref{l2.5}.

\medskip

\begin{lemma}\label{G:g6}
There exists $c=c(D,\alpha)>0$ such that for every $x,y \in D$ with
$r(x,y) < \eps_1$,
$$
g(A_{r(x,y)}(Q_x)) \wedge  g(A_{r(x,y)}(Q_y)) \,\ge\, c\,
r(x,y)^\gamma.
$$
\end{lemma}
\pf By symmetry, we show for $A= A_{r(x,y)}(Q_x)$ only. Note that
$g(\cdot)=G_D(\cdot,z_0)$ is harmonic in $D\cap B(Q_x, 2\eps_1)$.
Since $r(x,y) < \eps_1$, by Lemma \ref{l2.5} (recall $\eps_1 =
R/(12M)$),
$$
g(A)=G_D(A,z_0) \,\ge\, c\,\left(\frac{r(x,y)}{ 2\eps_1 }
\right)^{\gamma}G_D(A_{2\eps_1}(Q_x),z_0).
$$
Note that $ \rho_D(z_0) \ge  {R}/{M}=12 \eps_1$ and $
\rho_D(A_{2\eps_1}(Q_x)) > {2 \eps_1}/{M}. $ Thus by Lemma
\ref{G:g3}, $ G_D(A_{2\eps_1}(Q_x),z_0) >c_1 > 0$. \qed

\medskip

\begin{lemma}\label{G:g7}
There exists $c=c(D,\alpha)>0$ such that for every $x,y, z \in D$
and $(A_{x,y}, A_{y,z}) \in \BB(x,y) \times \BB(y,z)$
$$
\frac{g(A_{y,z})}{g(A_{x,y})}  \,\le\, c\,
\left(\frac{r(y,z)^\gamma}{r(x,y)^\gamma} \vee 1 \right).
$$
\end{lemma}
\pf Note that if $r(x,y) \ge \eps_1$, $g(A_{y,z}) \le C_1
=g(A_{x,y})$. So three cases will be considered:

a)  $r(x,y) < \eps_1$ and $  r(y,z) \ge \eps_1$: By Lemma
\ref{G:g6}, we have
$$
\frac{g(A_{y,z})}{g(A_{r(x,y)}(Q_y))}  \,\le\, c\,
\frac{C_1}{r(x,y)^\gamma}  \,\le\, c\, C_1 \eps_1^{-\gamma}
\frac{r(y,z)^\gamma}{r(x,y)^\gamma}.
$$

b) $r(y,z)\le r(x,y) < \eps_1$: Then $A_{r(y,z)}(Q_y) \in D \cap
B(Q_y, r(x,y))$. Thus by Lemma \ref{C:c_L}, $g(A_{r(y,z)}(Q_y))
\le c g(A_{r(x,y)}(Q_y)))$.

c)  $r(x,y) < r(y,z) < \eps_1$: By Lemma \ref{l2.5},
$$
\frac{g(A_{r(y,z)}(Q_y))}{g(A_{r(x,y)}(Q_y))}  \,\le\, c \,
\frac{r(y,z)^\gamma}{r(x,y)^\gamma} .
$$

In all three cases, we  apply Lemma \ref{G:g1} and we have proved
the lemma.  \qed

\medskip

\begin{lemma}\label{l:II}
There exists a constant $c=c(D,\alpha) >0$ such that for every  $x,
y, z, w\in D$ and $(A_{x,y}, A_{z,w}, A_{x,w}) \in \BB(x,y) \times
\BB(z,w) \times \BB(x,w)$,
$$
\frac{g(y)g(z) g(A_{x,w})^2}{g(A_{x,y})^2 g(A_{z,w})^2} \,\le\,
c\,\left(\frac{r(x,w)}{r(x,y) } \vee 1 \right)^\gamma
\left(\frac{r(x,w)}{r(z,w)}\vee 1 \right)^\gamma.
$$
\end{lemma}
\pf Note that
$$
\frac{g(y)g(z) g(A_{x,w})^2}{g(A_{x,y})^2 g(A_{z,w})^2} \,=\,
\frac{g(y)g(z)}{ g(A_{x,y}) g(A_{z,w})}
\left(\frac{g(A_{x,w})}{ g(A_{x,y})}\right)
\left(\frac{g(A_{x,w})}{ g(A_{z,w})}\right).
$$
By Lemma \ref{G:g5}, we have
$$
\frac{g(y)g(z)}{ g(A_{x,y}) g(A_{z,w})} \, \le \, c.
$$
On the other hand, applying Lemma \ref{G:g7} gives
$$
\left(\frac{g(A_{x,w})}{ g(A_{x,y})}\right)
\left(\frac{g(A_{x,w})}{ g(A_{z,w})}\right) \,\le\,
c\,\left(\frac{r(x,w)}{r(x,y) } \vee 1 \right)^\gamma
\left(\frac{r(x,w)}{r(z,w)}\vee 1 \right)^\gamma.
$$
\qed

\medskip

\begin{lemma}\label{l:in}
For every $a, b, c >0$, we have
$$
\left(\frac{a}{b} \vee 1\right)+ \left(\frac{a}{c} \vee 1\right)
\,\le\, 2 \left(\frac{a}{b} \vee 1 \right) \left(\frac{a}{c} \vee
1\right).
$$
\end{lemma}
\pf By considering each case, it is clear. So we omit details.
\qed

\medskip

\begin{lemma}\label{l:I}
There exists a constant $c=c(D,\alpha) >0$ such that for every $x, y, z, w\in
D$ and $(A_{x,y}, A_{z,w}, A_{x,w}) \in \BB(x,y) \times \BB(z,w)
\times \BB(x,w)$,
$$
\frac{g(y)g(z) g(A_{x,w})^2}{g(A_{x,y})^2 g(A_{z,w})^2} \,\le\,
c\,\left(\frac{r(y,z)}{r(x,y) } \vee 1 \right)^\gamma
\left(\frac{r(y,z)}{r(z,w)}\vee 1 \right)^\gamma.
$$
\end{lemma}
\pf From Lemma \ref{g:tri3}, we get
\begin{eqnarray}
\frac{g(y)g(z) g(A_{x,w})^2}{g(A_{x,y})^2 g(A_{z,w})^2} &\le&
c\,\frac{g(y)g(z)}{g(A_{x,y})^2 g(A_{z,w})^2}(  g(A_{x,y})^2 +
g(A_{y,z})^2 +
g(A_{z,w})^2 )\nonumber\\
&=& c\left(\frac{g(y)g(z)}{g(A_{z,w})^2} +\frac{g(y)g(z)}{g(A_{x,y})^2}
+\frac{g(y)g(z)g(A_{y,z})^2} {g(A_{x,y})^2 g(A_{z,w})^2}\right).
\label{e:3G33}
\end{eqnarray}
By applying Lemma \ref{G:g5} to both $y$ and $z$, we have that
$g(y)  \le c g(A_{x,y})  $ and $g(z)  \le c g(A_{z,w})$. Thus
(\ref{e:3G33}) is less than or equal to
\begin{eqnarray*}
&&c\frac{g(y)}{g(A_{z,w})} +c\frac{g(z) }{g(A_{x,y})}
+c^2\left(\frac{g(A_{y,z})} {   g(A_{x,y})}\right)
\left(\frac{g(A_{y,z})} {   g(A_{z,w})}\right)\\
&&\le c\frac{g(y)}{g(A_{z,w})} +c\frac{g(z) }{g(A_{x,y})} +c_1
\left(\frac{r(y,z)^\gamma}{r(x,y)^\gamma} \vee 1 \right)
\left(\frac{r(y,z)^\gamma}{r(z,w)^\gamma} \vee 1 \right).
\end{eqnarray*}
We have used Lemma \ref{G:g7} in the last inequality above.
Moreover, by Lemma \ref{G:g5} and Lemma \ref{G:g7},
$$
\frac{g(y)}{g(A_{z,w})} =\left( \frac{g(y)}  {g(A_{y,z})}
\right)\left( \frac{g(A_{y,z})}  {g(A_{z,w})}
 \right)\, \le\, c \,\left(\frac{r(y,z)^\gamma}{r(z,w)^\gamma} \vee 1 \right)
$$
and
$$
\frac{g(z)}{g(A_{x,y})} =\left( \frac{g(z)}  {g(A_{y,z})}
\right)\left( \frac{g(A_{y,z})}  {g(A_{x,y})}
 \right)\, \le\, c \, \left(\frac{r(y,z)^\gamma}{r(x,y)^\gamma} \vee 1 \right).
$$
We conclude that (\ref{e:3G33}) is less than or equal to
$$
c_2\, \left(\left(\frac{r(y,z)^\gamma}{r(x,y)^\gamma} \vee 1
\right)+ \left(\frac{r(y,z)^\gamma}{r(z,w)^\gamma} \vee 1 \right)
+ \left(\frac{r(y,z)^\gamma}{r(x,y)^\gamma} \vee 1 \right)
\left(\frac{r(y,z)^\gamma}{r(z,w)^\gamma} \vee 1 \right)\right)
$$
 for some $c_2=c_2(D, \alpha)>0$.
Using Lemma \ref{l:in}, we have finished the proof. \qed

\medskip

Now we are ready to prove the main result of this paper.

\medskip

\noindent{\bf Proof} of \noindent{\bf Theorem} \ref{t:3Gt}.
 Let
$$
G(x,y,z,w)\,:=\,\frac{G_D(x,y) G_{D}(z,w)} { G_{D}(x,w)} \quad
\mbox{and} \quad H(x,y,z,w)\,:=\,\frac{ |x-w|^{d-\alpha}}
{|x-y|^{d-\alpha} |z-w|^{d-\alpha}   }.
$$

\begin{enumerate}
\item If  $|x-w| \le  \rho_D(x) \wedge \rho_D(w)$, by Lemma
\ref{G:g3}, $ G_D(x,w)  \ge c |x-w|^{-d+\alpha}. $ Thus we have $
G(x,y,z,w) \le cH(x,y,z,w). $

\item If $|y-z| \leq \rho_D(y) \wedge \rho_D(z)$, then by Lemma
\ref{G:g3} $ G_D(y,z)  \ge c |y-z|^{-d+\alpha}. $ Moreover, from
the classical 3G theorem (Theorem 6.1 in \cite{SW}), we obtain
that there exists a constant $c=c(D, \alpha) >0$ such that
\begin{eqnarray}
&& G(x,y,z,w) \,=\,
\frac{G_D(x,y) G_{D}(y,z)} { G_{D}(x,z)}  \frac{G_D(x,z) G_{D}(z,w)} { G_{D}(x,w)}\frac{1}{ G_{D}(y,z) }\nonumber\\
&&\le \,c \,\frac{|x-z|^{d-\alpha}} {|x-y|^{d-\alpha}
|y-z|^{d-\alpha}}\frac{|x-w|^{d-\alpha}}
{|x-z|^{d-\alpha} |z-w|^{d-\alpha}}\frac{1}{   G_{D}(y,z)   }\nonumber\\
&&=\, c \, \frac{ |x-w|^{d-\alpha}} {|x-y|^{d-\alpha}
|y-z|^{d-\alpha}|z-w|^{d-\alpha} }\frac{1}{ G_{D}(y,z)}.
\label{3G_est53}
\end{eqnarray}
Thus, we have $ G(x,y,z,w) \le cH(x,y,z,w). $
    \item
    Now we assume that $|x-w| >  \rho_D(x) \wedge \rho_D(w)$ and $|y-z| > \rho_D(y) \wedge
    \rho_D(z)$. Since $\rho_D(x) \vee \rho_D(w) \le \rho_D(x) \wedge \rho_D(w) +|x-w|$, using the assumption
    $ \rho_D(x) \wedge \rho_D(w)< |x-w|$, we obtain $r(x,w) <2|x-w|$. Similarly, $r(y,z)<2|y-z|$.

Let $A_{x,w} \in \BB(x,w),\ A_{x,y} \in
\BB(x,y)$ and $ A_{z,w} \in \BB(z,w)$. Applying Lemmas \ref{l:II} and
\ref{l:I} to (\ref{3G_est}), we have
\begin{align}\label{G:almost}
&G(x,y,z,w)\, \le\, c\, \frac{g(y)g(z) g(A_{x,w})^2}{g(A_{x,y})^2
g(A_{z,w})^2} H(x,y,z,w) \notag\\
&\le\, c \left(\frac{r(y,z) \wedge r(x,w)}{r(x,y)} \vee 1
\right)^\gamma \left(\frac{r(y,z)\wedge r(x,w)}{r(z,w)} \vee 1
\right)^\gamma H(x,y,z,w) .
\end{align}
Now applying the fact that $r(x,w)<2|x-w|$,
$r(y,z)<2|y-z|$, $r(x,y) \ge |x-y|$ and $r(z,w) \ge
|z-w|$,  we have (\ref{3G}).
\end{enumerate}
We have proved the theorem. \qed

\medskip

\begin{remark}\label{rmk1}
{\rm By letting $y=z$ in (\ref{3G}), we recover the classical 3G
theorem, i.e.,
$$
\frac{G_D(x,z) G_{D}(z,w)} { G_{D}(x,w)} \, \le\, \,c\,
\frac{|x-w|^{d-\alpha}} {|x-z|^{d-\alpha} |z-w|^{d-\alpha}}
\,\le\, \,c_1\left( \frac{1} {|x-z|^{d-\alpha}} +\frac{1}
{|z-w|^{d-\alpha}} \right).
$$
}
\end{remark}

\medskip

\begin{remark}\label{rmk0}
{\rm See (3.2)-(3.3) in \cite{CS7} for a generalized version of
3G-estimate on bounded $C^{1,1}$ domains.
Using Green function estimates in \cite{CS1}, instead of Lemma
\ref{l2.5},
it is easy to see
that $\gamma=\alpha/2$ in (\ref{3G})
for  bounded $C^{1,1}$ domains.
}
\end{remark}

\medskip

\begin{remark}\label{rmk2}
{\rm The inequality in (\ref{3G}) is not true even in smooth
domains without the first two factors in the right-hand side of
(\ref{3G}); If $D$ is smooth, by Green function estimate (for
example, see \cite{CS1}), there exists a positive constant
$c=c(D,\alpha)$ such that
\begin{equation}\label{e:rmk2}
c \left(\frac{\rho_D(y)\rho_D(z)|x-w|^2}
{r(x,y)^2r(z,w)^2}\right)^{\frac{\alpha}2} \frac{
|x-w|^{d-\alpha}} {|x-y|^{d-\alpha} |z-w|^{d-\alpha}   }   \,\le\,
\frac{G_D(x,y) G_{D}(z,w)} { G_{D}(x,w)}.
\end{equation}
Suppose the inequality in (\ref{3G}) is true  without the factors.
Then (\ref{e:rmk2}) implies that there exists a positive constant
$c_1=c_1(D,\alpha)$ such that for distinct $x,y,z,w$ in $D$
\begin{equation}\label{e:rmk3}
\frac{\rho_D(y)\rho_D(z)|x-w|^2} {r(x,y)^2r(z,w)^2} \,\le\, c_1.
\end{equation}
But (\ref{e:rmk3}) can not be true. In fact, for points $x \not=w
\in D$ near $\partial D$, choose $ y\not=z$ such that $x \neq y$, $z
\neq w$, $r(x,y)= \rho_D(x) =\rho_D(y) \ge |x-y|$ and $r(z,w)=
\rho_D(z)= \rho_D(w) \ge |z-w|$. Then,
$$
\frac{\rho_D(y)\rho_D(z)|x-w|^2} {r(x,y)^2r(z,w)^2} \,=\,
\frac{|x-w|^2} {\rho_D(x)\rho_D(w)}\, \to \infty \quad \mbox{ as
}\, x \to \partial D.
$$
}
\end{remark}

\section{Concrete sufficient conditions for the
non-local Kato class and relativistic stable
process}\label{sec:kato} Throughout this section $D$ is a bounded
$\kappa$-fat open set. In this section, we recall the new
(non-local) Kato class introduced by Chen and Song
(\cite{C,CS6,CS7}). Then we apply the generalized 3G theorem to
establish some concrete sufficient condition for functions to be in
the non-local Kato class. Using this sufficient condition, we
discuss some properties of relativistic stable process in bounded
$\kappa$-fat open sets.

We will use the definitions for the new Kato class not only for
symmetric stable process, but also for other symmetric Hunt
processes, which will be specified later. So we will denote $Y$ for
those processes in $D$,  $G(x,y)$ for the Green function for $Y$ and
will state the following definition for $Y$. Note that the processes
in this paper are irreducible transient symmetric Hunt processes
satisfying the assumption at the beginning of section 3.2 in
\cite{C}.

We call a positive measure $\mu$ on $D$ a smooth measure of $Y$ if
there is a positive continuous additive functional (PCAF in
abbreviation) $A$ of $Y$ such that
\begin{equation}\label{eqn:Revuz1}
\int_D f(x) \mu (dx) = \uparrow \lim_{t\downarrow 0} \int_D \E_x
\left[ \frac1t \int_0^t f(Y_s) dA_s \right] dx
\end{equation}
for any Borel measurable function $f\geq 0$. Here $\uparrow
\lim_{t\downarrow 0}$ means the quantity is increasing as
$t\downarrow 0$. The measure $\mu$ is called the Revuz measure of
$A$. For a signed measure $\mu$, we use $\mu^+$ and $\mu^-$ to
denote its positive and negative parts respectively. If $\,\mu^+$and
$\mu^{-}$ are smooth measures of $Y$ and $A^{+}$ and $A^{-}$ are
their corresponding PCAFs of $Y$, then we say the continuous
additive functional $A:=A^+-A^{-}$ of $Y$ has (signed) Revuz measure
$\mu$.

\medskip

\begin{defn}\label{df:5.1}
Suppose that $A$ is a continuous additive functional of $Y$ with
Revuz measure $\nu$. Let $A^{+}$ and $A^{-}$ be the PCAFs (positive
continuous additive functionals) of $Y$ with Revuz measures
$\nu^{+}$ and $\nu^{-}$ respectively. Let $|A|= A^{+}+A^{-}$ and
$|\nu|=\nu^{+}+\nu^{-}$.

\begin{description}
\item{(1)} The measure $\nu$ (or the continuous additive
functional $A$) is said to be in the class $\S_\infty (Y)$ if for
any $\eps>0$ there are a Borel subset $K=K(\eps)$ of finite $|
\nu|$-measure and a constant $\delta = \delta (\eps) >0 $ such that
$$ \sup_{(x, z)\in (D\times D) \setminus \{x = z\}} \int_{D\setminus K}
\frac{G(x, y)G(y, z) }{ G(x, z)} \, |\nu | (dy)\le\eps $$ and that,
for all measurable set $B\subset K$ with $| \nu |(B)<\delta$,
$$ \sup_{(x, z)\in (D\times D) \setminus \{x =z\}} \int_B \frac{G(x,
y)G(y, z) }{ G(x, z)} \, |\nu |  (dy) \le\eps. $$

\item{(2)} A function $q$ is said to be in the class $\S_\infty
(Y)$, if $\nu(dx) :=q(x)dx$ is in the class $\S_\infty (Y)$.
\end{description}
\end{defn}
\medskip
Let $(N, H)$ be a L\'evy system for $Y$ in $D$, that is, for every
$x\in D$, $N(x, dy)$  is a kernel on $(D_\partial, {\cal
B}(D_\partial))$, where $\partial $ is the  cemetery point for
process $Y$ and $D_\partial = D\cup \{ \partial \}$, and $H_t$ is a
positive continuous additive functional of $Y$ with bounded
1-potential such that  for any nonnegative Borel function $f$ on $D
\times D_\partial $ that vanishes along the diagonal $d$,
$$
\E_x\left(\sum_{s\le t}f(Y_{s-}, Y_s) \right) =
\E_x\left(\int^t_0\int_{D_\partial} f(Y_s, y) N(Y_s, dy)dH_s\right)
$$
for every $x\in D$ (see \cite{Sharpe} for details). We let
$\mu_{H}(dx)$ be the Revuz measure for $H$.
\medskip

\begin{defn}\label{df:5.2}
Suppose $F$ is a bounded function on $D\times D$ vanishing on the
diagonal. Let $$ \mu_{|F|} (dx):=\left( \int_D |F(x, y)| N(x,
dy)\right) \mu_{H}(dx). $$

(1) $F$ is said to be in the class $\A_\infty(Y)$ if for any
$\eps>0$ there are a Borel subset $K=K(\eps)$ of finite
$\mu_{|F|}$-measure and a constant $\delta = \delta (\eps) >0 $ such
that $$ \sup_{(x, w)\in (D\times D) \setminus \{x =w\}}
\int_{(D\times D) \setminus (K\times K)} G(x, y) \frac{|F(y, z)|G(z,
w) }{ G(x, w)}N(y,dz)\mu_{H}(dy) \,\le\,\eps $$ and that, for all
measurable sets $B\subset K$ with $\mu_{|F|} (B)<\delta$, $$
\sup_{(x, w)\in (D\times D) \setminus \{x =w\}} \int_{(B\times
D)\cup (D\times B)} G(x, y) \frac{|F(y, z)| G(z, w) }{G(x,
w)}N(y,dz)\mu_{H}(dy)\,\le\,\eps. $$

(2) $F$ is said to be in the class $\A_2 (Y)$ if $F\in \A_\infty(Y)$
and if the measure  $\mu_{|F|}$ is in $\S_\infty (Y)$.

\end{defn}

\medskip

By (\ref{KDF}),  L\'{e}vy system $(N^D, H^D)$ for symmetric stable
process $X^D$ can be chosen to be
\begin{equation}\label{e:R}
N^D(x,dy)=2J(x,y)dy=\frac{{\cal A}(d,-\alpha)}{|x-y|^{d+\alpha}}dy
\quad \mbox{in} \,\,D  \quad\quad \hbox{and} \quad\quad  H^D_t=t,
\end{equation}
where $J(x,y)$ is given in (\ref{J}). By (\ref{eqn:Revuz1}), the
Revuz measure $\mu_{H^D} (dx)$ for $H^D$ is simply the Lebesgue
measure $dx$ on $D$. (see Theorem 4.5.2, Theorem 4.6.1, Lemmas
5.3.1--5.3.3 and Theorem A.3.21 in \cite{FOT} for the relation
between (regular) Dirichlet form and L\'{e}vy system.)

Using the generalized 3G theorem (Theorem \ref{t:3Gt}), now we can
give a concrete condition for a function $F$ to be in class
$\A_2(X^D)$ on any bounded $\kappa$-fat open set $D$. For $w \in
D$, We denote by $\E_x^w$ the expectation for the conditional
process obtained from  $X^D$ through Doob's $h$-transform with
$h(\cdot)=G_D(\cdot, w)$ starting from $x\in D$.

\medskip

\begin{thm}\label{t:A_2}
If $D$ is a bounded $\kappa$-fat open set and $F$ is a function on
$D \times D$ with
\beq\label{F}
|F(x,y)| \leq c|x-y|^{\beta}
\eeq
for some $\beta > \alpha$ and $c
> 0 $, then $F \in \A_2(X^D)$ and
\beq\label{bound}
 \sup_{(x,w)\in D \times D} \E^w_x \left[ \sum_{s <
\tau_D}F(X^D_{s-}, X^D_s) \right] \, < \, \infty.
\eeq
\end{thm}
\pf First let us show that $F \in \A_{\infty}$. Since $D$ is
bounded, by the Generalized 3G theorem (Theorem \ref{t:3Gt}), there
exists a positive constant $c=c(D,\alpha)$ such that
\begin{eqnarray}
&& \frac{G_D(x,y) G_{D}(z,w)} { G_{D}(x,w)} \nonumber\\
&& \le \,c\, \left( \frac{1}{|x-y|^{d-\alpha+\gamma}}+
\frac{1}{|z-w|^{d-\alpha+\gamma}} + \frac{|y-z|^{d-\alpha}}
{|x-y|^{d-\alpha} |z-w|^{d-\alpha}} +\frac{|y-z|^{d-\alpha+\gamma}}
{|x-y|^{d-\alpha}
|z-w|^{d-\alpha+\gamma}}\right.\nonumber\\
&&~~~~~~~~\left.+ \frac{|y-z|^{d-\alpha+\gamma}}
{|x-y|^{d-\alpha+\gamma} |z-w|^{d-\alpha}}+
\frac{|y-z|^{d-\alpha+2\gamma}} {|x-y|^{d-\alpha+\gamma}
|z-w|^{d-\alpha+\gamma}}\right).\label{3G_est6}
\end{eqnarray}
Thus, by the assumption on $F$, (\ref{F}), and (\ref{e:R}), we
obtain
\begin{eqnarray}
&&\frac {G_D(x,y)|F(y,z)|  G_D(z,w)}{G_D(x,w)} \frac{N^D(y,dz)}{dz}
\,\le\,c\,\frac {G_D(x,y)
G_D(z,w)}{G_D(x,w)}\frac{1}{|y-z|^{d+\alpha-\beta}}
\nonumber\\
&&\le \,c\, \left(
\frac{1}{|x-y|^{d-\alpha+\gamma}|y-z|^{d+\alpha-\beta}}
+\frac{1}{|y-z|^{d+\alpha-\beta}|z-w|^{d-\alpha+\gamma}}
\right.\nonumber\\
&&~~~~ + \frac{1} {|x-y|^{d-\alpha}|y-z|^{2\alpha - \beta }
|z-w|^{d-\alpha}}+ \frac{1} {|x-y|^{d-\alpha}|y-z|^{2\alpha - \beta
-\gamma} |z-w|^{d-\alpha+\gamma}}
\nonumber\\
&&~~~~\left. + \frac{1} {|x-y|^{d-\alpha+\gamma} |y-z|^{2\alpha -
\beta-\gamma}|z-w|^{d-\alpha}}+ \frac{1}
{|x-y|^{d-\alpha+\gamma}|y-z|^{2\alpha - \beta-2\gamma}
|z-w|^{d-\alpha+\gamma}}\right).\label{eqn:4.10}
\end{eqnarray}

Since $\gamma<\alpha<\beta$,
\begin{eqnarray*}
&& \{ (y,z)\mapsto |x-y|^{\alpha-\gamma-d} |y-z|^{\beta-\alpha-d}, \, x\in D\}, \\
&&\{  (y,z)\mapsto |z-w|^{\alpha-\gamma-d} |y-z|^{\beta-\alpha-d},\,
w \in D\}
\end{eqnarray*}
are uniformly integrable over cylindrical sets of the form $B\times
D$ and $D\times B$, for any Borel set $B \subset D$.

Now let us show that the following family of functions are uniformly
integrable over cylindrical sets of the form $B\times D$ and
$D\times B$:
\begin{eqnarray}
&&\{ (y,z)\mapsto |x-y|^{\alpha-d} |z-w|^{\alpha-d}
|y-z|^{\beta-2\alpha}, \,
\, x, w \in D\},\label{f1}\\
&&\{ (y,z)\mapsto |x-y|^{\alpha-d} |z-w|^{\alpha-\gamma-d}
|y-z|^{\beta-2\alpha+\gamma}, \,
\, x, w \in D\},\label{f2}\\
&&\{ (y,z)\mapsto |x-y|^{\alpha-\gamma-d} |z-w|^{\alpha-d}
|y-z|^{\beta-2\alpha+\gamma}, \,
\, x, w \in D\},\label{f3}\\
&&\{ (y,z)\mapsto |x-y|^{\alpha-\gamma-d} |z-w|^{\alpha-\gamma-d}
|y-z|^{\beta-2\alpha+2\gamma}, \, \, x, w \in D\}.\label{f4}
\end{eqnarray}
We apply Young's inequality $ab \le \frac{1}{p}a^p+\frac{1}{q}b^q$
for some $p, q >1$ satisfying
 \beq \label{pq}
 p^{-1}+q^{-1}=1.
 \eeq
\begin{enumerate}
\item
 First, let us consider the function (\ref{f1}) when the exponent of $|y-z|$ is negative, i.e., $\beta
< 2\alpha$. Otherwise it is obvious that the function (\ref{f1}) is
uniformly integrable since $D$ is bounded. Applying Young's
Inequality, we obtain
\begin{align*}
& |x-y|^{\alpha-d}|z-w|^{\alpha-d}|y-z|^{\beta-2\alpha}=
 \bigl(|x-y||z-w|\bigr)^{\alpha-d}|y-z|^{\beta-2\alpha}\\
& \leq \frac{1}{p}\Bigl(|x-y||z-w|\Bigr)^{(\alpha
-d)p}+\frac{1}{q}|y-z|^{(\beta -2\alpha)q}.
\end{align*}
It suffices to find $p, q >1$ satisfying (\ref{pq}) and
\beq \label{pq-f1} (d-\alpha)p<d,\quad (2\alpha -\beta )q<d. \eeq
Choosing $p$ in the interval
$$
\left( 1 \vee \frac{d}{d-2\alpha +\beta } \right) < p <
\frac{d}{d-\alpha},
$$
we can finish this case. Note that this interval is not empty since
$\dfrac{d}{d-2\alpha +\beta } <\dfrac{d}{d-\alpha }$ by $\beta
\wedge d> \alpha$ and $\dfrac{d}{d-\alpha }>1$.

\item For the function (\ref{f4}), we replace $\alpha$ by $\alpha
-\gamma $ in the previous case. Here note that $0<\alpha -\gamma
<\beta \wedge d$.

\item Now let us consider the function (\ref{f2}) only when the
exponent of $|y-z|$ is negative, i.e., $\beta +\gamma < 2\alpha$.
\begin{align*}
&|x-y|^{\alpha-d}
|z-w|^{\alpha-\gamma-d}|y-z|^{\beta-2\alpha+\gamma}\\
&=\left(|x-y||z-w|\right)^{\alpha-d}\left(|z-w|^{-\gamma}|y-z|^{\beta-2\alpha+\gamma}\right)\\
&\le
\frac{1}{p}\Bigl(|x-y||z-w|\Bigr)^{(\alpha-d)p}+\frac{1}{q}|z-w|^{-\gamma
q}|y-z|^{(\beta-2\alpha+\gamma)q}.
\end{align*}
We find $p, q >1$ satisfying (\ref{pq}) and
\beq \label{pq-f2} (d-\alpha)p<d,\quad \gamma q<d,\quad (2\alpha
-\beta-\gamma)q<d. \eeq
Choose $p$ in the interval
$$
\left( 1 \vee \frac{d}{d-\gamma} \vee \frac{d}{d-2\alpha
+\beta+\gamma} \right) < p < \frac{d}{d-\alpha} .
$$
This interval is not empty. Clearly, $\dfrac{d}{d-\alpha}>1$. Also,
$\dfrac{d}{d-\alpha}>\dfrac{d}{d-\gamma}$ since $0<\gamma
<\alpha<d$. Finally, $\dfrac{d}{d-\alpha}>\dfrac{d}{d-2\alpha
+\beta+\gamma}$ since
$d>\alpha>\alpha-(\beta-\alpha)-\gamma=2\alpha-\beta-\gamma$,
$\alpha<\beta$.

\item Analogously to the previous case, the function (\ref{f3})
can be considered.

\end{enumerate}

Thus $F$ is in $\A_\infty(X^D)$ and, moreover,
$$ \sup_{(x,w)\in D
\times D} \int_D \int_D \frac{G_D(x,y) G_{D}(z,w)} { G_{D}(x,w)}
F(y,z)|y-z|^{-d-\alpha}dydz \, < \, \infty.
$$
Note that, by Proposition 3.3 in \cite{CS7} and
(\ref{e:R}),
$$ \E_x^w\left[ \sum_{s < \tau_D}F(X^D_{s-}, X^D_s) \right]={\cal A}
(d,-\alpha) \int_D \int_D \frac{G_D(x,y) G_{D}(z,w)} { G_{D}(x,w)}
F(y,z)|y-z|^{-d-\alpha}dydz. $$
Therefore we have (\ref{bound}).

Since
$$ \mu_{|F|} (dx) = {\cal A}(d, -\alpha)
\left( \int_D |F(x,y)| |x-y|^{-\alpha-d}dy \right) dx \leq  \left(
\int_D  c\, |x-y|^{\beta-\alpha-d}dy \right) dx \leq c \, dx,
$$
it follows from (\ref{3G_est6}) with $y=z$ that $\mu_{|F|}\in
\S_\infty(X^D)$ and therefore  $F$ is in $\A_2(X^D)$. \qed

\medskip

Fix $x_0 \in D$ and set
$$M_D(x,y) := \frac {G_D(x,y)}{G_D(x_0 , y)} ~,~~~~ x,y \in D.$$
It is shown in \cite{SW} that
\begin{equation}\label{e:Mcon}
M_D(x,z):= \lim_{y\rightarrow z \in \partial D} M_D(x,y) \quad
\mbox{exists for every } z\in \partial D,
\end{equation}
which is called the Martin kernel of $D$, and that $M_D(x, z)$ is
jointly  continuous in $D\times \partial D$. For each $z \in
\partial D$, set $M_D(x,z)=0$ for $x \in D^c$. For each $z\in \partial D$,
$x\mapsto M_D(x, z)$ is a singular harmonic function for $X$ in $D$
(see pages 415-416 of \cite{SW} for details). We get the following
as a corollary of the generalized 3G-inequality.

\medskip

\begin{corollary}\label{c:GMM}
For every $x,y,z \in D$ and $w \in \partial D$, there exists a
constant $c=c(D,\alpha) >0$ such that
\begin{equation} \label{e:GMM}
\frac{G_D(x,y) M_{D}(z,w)} { M_{D}(x,w)} \le c
\left(\frac{|x-w|\wedge |y-z|}{|x-y|  } \vee 1 \right)^\gamma
\left(\frac{|x-w|\wedge |y-z|}{|z-w|} \vee 1 \right)^\gamma
\frac{|x-w|^{d-\alpha}} {|x-y|^{d-\alpha} |z-w|^{d-\alpha}}.
\end{equation}
\end{corollary}

\medskip

Let $H^+(X^D)$ be the collection of all positive singular harmonic
function of $X$ in $D$. We denote $\E_x^h$ the expectation for the
conditional process obtained from  $X^D$ through Doob's
$h$-transform with $h(\cdot)$ starting from $x\in D$.
i.e.,
$$
\E_x^h \left[f(X^D_t) \right]\,=\,
\E_x\left[\frac{h(X^D_t)}{h(x)}f(X^D_t) \right]
.$$
When
$h(\cdot)= M_D(\cdot, w)$ with $w \in \partial D$, we use $\E_x^w$
instead of $\E_x^{M_D(\cdot, w)}$.

\medskip

The next result can be proved easily with an argument similar to the
argument in the proof of Proposition 3.3 in \cite{CS6}. We put the
details for the reader's convenience.
\begin{prop}\label{c:conh}
For every $\beta > \alpha$,
\begin{equation} \label{e:conh}
\sup_{x\in D, h \in H^+(X^D)} \E^h_x \left[ \sum_{s <
\tau_D}|X^{D}_{s-}-X^{D}_s|^{\beta} \right] \, < \, \infty.
\end{equation}
\end{prop}
\pf For every $h \in H^+(X^D)$, $M_t:=h(X^D_s)/h(X^D_0)$ is a
supermartingale multiplicative functional of $X^D$. It follows from
Section 62 of \cite{Sharpe}
 that
\begin{eqnarray*}
&&\E^h_x \left[ \sum_{s \le t}|X^{D}_{s-}-X^{D}_s|^{\beta} \right]\\
&&= \E_x \left[ \sum_{s \le t}|X^{D}_{s-}-X^{D}_s|^{\beta}
\left(\int_0^{\tau_D} 1_{\{s<r \le \tau_D
\}}d(-M_r)+M_{\tau_D}\right);\, t < \tau_D
\right]\\
&&= \E_x \left[ \sum_{s \le t}|X^{D}_{s-}-X^{D}_s|^{\beta} M_s;\, t
< \tau_D
\right]\\
&&= \frac1{h(x)}\E_x \left[ \sum_{s \le
t}|X^{D}_{s-}-X^{D}_s|^{\beta} h(X^D_s);\, t < \tau_D \right]\\&&=
\frac1{h(x)}\E_x \left[ \int_0^t \int_D |X^{D}_{s}-y|^{\beta} h(y)
N^D( X^{D}_{s}, dy) d H_s^D \right].
\end{eqnarray*}
Thus by (\ref{e:R})
\begin{eqnarray}
\E^h_x \left[ \sum_{s < \tau_D}|X^{D}_{s-}-X^{D}_s|^{\beta} \right]
&=& \frac1{h(x)}\E_x \left[ \int_0^{\tau_D} \int_D
|X^{D}_{s}-y|^{\beta} h(y) N^D( X^{D}_{s}, dy) d H_s^D
\right]\nonumber\\
 &=& {\cal A} (d, \, -\alpha)
\int_D \int_D \frac{G_D(x,y) h(z)} { h(x)}
|y-z|^{-d-\alpha+\beta}dydz.\label{e:pr}
\end{eqnarray}
Note that since, by Corollary \ref{c:GMM},
$$
\frac {G_D(x,y) M_D(z,w)}{M_D(x,w)}|y-z|^{-d-\alpha+\beta}$$ is
bounded above by (\ref{eqn:4.10}), following the argument in the
proof of Theorem \ref{t:A_2}, we get
\begin{equation} \label{e:GMM2}
\sup_{(x,w)\in D \times \partial D} \int_D \int_D \frac{G_D(x,y)
M_{D}(z,w)} { M_{D}(x,w)}  |y-z|^{-d-\alpha+\beta}dydz \, < \,
\infty.
\end{equation}
Thus (\ref{e:GMM2}) and (\ref{e:pr}) with $h(\cdot)= M_D(\cdot, w)$
imply that
\begin{equation} \label{e:GMM3}
\sup_{(x,w)\in D \times \partial D} \E^w_x \left[ \sum_{s <
\tau_D}|X^{D}_{s-}-X^{D}_s|^{\beta} \right] \, < \, \infty.
\end{equation}

On the other hand, from the Martin representation ((4.1) in
\cite{SW}), for every $h \in H^+(X^D)$, there exists a finite
measure $\mu$ on $\partial D$ such that $ h(x) =\int_{\partial D}
M_D(x,w) \mu(dw).$ Thus by (\ref{e:pr}) and (\ref{e:GMM3}),
\begin{eqnarray*}
h(x) \E^h_x \left[ \sum_{s < \tau_D}|X^{D}_{s-}-X^{D}_s|^{\beta}
\right]
&=&\int_{\partial D} M_D(x,w) \E^w_x \left[ \sum_{s <  \tau_D}|X^{D}_{s-}-X^{D}_s|^{\beta} \right]  \mu(dw)\\
&\le& c \,\int_{\partial D} M_D(x,w) \mu(dw)\, = \,c\, h(x).
\end{eqnarray*}
Since $0<h(x)< \infty$, we have (\ref{e:conh}).\qed

\medskip

For the remainder of this section, we will discuss some properties
of relativistic $\alpha$-stable process in bounded $\kappa$-fat open
sets. This process have been studied in \cite{CK2, CS4, K, R}. We
will recall the definition of  relativistic $\alpha$-stable process
and its basic  properties. We will spell out some of the details for
the readers, who are unfamiliar with the non-local Feynman-Kac
perturbations (cf. \cite{C,CS7, CK2}).

For $0<\alpha<2$, a relativistic $\alpha$-stable process $X^m$  in
$\R^d$ is a L\'evy process whose characteristic function is given by
$$
\E \left[\exp\left(i \xi \cdot (X^m_t-X^m_0) \right) \right] =
\exp\left(-t \left( (|\xi|^2+ m^{2/\alpha} )^{\alpha/2}-m\right)
\right), \quad \xi \in \R^d,
$$
where $m>0$ is a constant. In other words, the relativistic
$\alpha$-stable process in $\R^d$ has infinitesimal generator $m-
(-\Delta +m^{2/\alpha} )^{\alpha /2}$. It is clear that $X^m$ is
symmetric with respect to the Lebesgue measure $dx$ on $\R^d$ and
that when the parameter $m$ degenerates to $0$, $X^m$ becomes a
symmetric $\alpha$-stable process $X$ on $\R^d$.

The Dirichlet form $({\cal Q}, D({\cal Q}))$ of $X^m$ is given by
$$
{\cal Q}(u, v):=\int_{\R^d} \hat{v} (\xi)\bar{\hat{u}}(\xi) \left(
(|\xi|^2+ m^{2/\alpha} )^{\alpha/2}-m\right) d\xi
$$
and
$$
D({\cal Q}):= \{ u\in L^2(\R^d):\int_{\R^d} |\hat{u} (\xi)|^2\left(
(|\xi|^2+ m^{2/\alpha} )^{\alpha/2}-m\right) d\xi <\infty\}.
$$
where $\hat{u} (\xi)=(2\pi)^{-d/2}\int_{\R^d}e^{i\xi\cdot y}u(y)dy$
(see Example 1.4.1 of \cite{FOT}). Recall $({\cal E}, {\cal F})$ from  (\ref{DF1}).
Similar to ${\cal E}_1$, we can
also define ${\cal Q}_1$. From (\ref{DF1}), we see that there exist
 positive
constants $c_1=c_1(m)$ and $c_2=c_2(m)$ such that
\begin{equation}\label{DFC}
c_1 {\cal E}_1(u, u)\le {\cal Q}_1(u, u)\le c_2 {\cal E}_1(u, u).
\end{equation}
Therefore $D({\cal Q})={\cal F}$. From now on we will use ${\cal F}$
instead of $D({\cal Q})$.

Like symmetric stable processes, relativistic $\alpha$-stable
processes can be obtained from Brownian motions through
subordinations. For the details, we refer our readers to \cite{R}.
By Lemma 2 in \cite{R}, the L\'{e}vy measure for relativistic
$\alpha$-stable process has the density
$$
\nu(x):=\frac{\alpha}{2 \Gamma (1-\frac{\alpha}2)} \int_0^{\infty}
(4\pi
u)^{-d/2}e^{-\frac{|x|^2}{4u}}e^{-m^{2/\alpha}}u^{-(1+\frac{\alpha}{2})}du.
$$
Using change of variables twice, first with $u=|x|^2v$ then with
$v=1/s$, we get
$$
\nu(x)=  {\cal A} (d, \, -\alpha) |x|^{-d-\alpha} \psi
(m^{1/\alpha}|x|)
$$
where
$$
\psi (r):= 2^{-(d+\alpha)} \, \Gamma \left( \frac{d+\alpha}{2}
\right)^{-1}\, \int_0^\infty s^{\frac{d+\alpha}{ 2}-1} e^{-\frac{s}{
4} -\frac{r^2}{ s} } \, ds,
$$
which is a smooth function of $r^2$ (see pages 276-277 of \cite{CS4}
for the details). Recall $J(x,y)$ from (\ref{J}). The L\'{e}vy
measure $\nu(x)dx$ determines the jumping measure $J^m$ of $X^m$:
$$
J^m(x,y):=\frac12 \nu(x-y) =\frac12 {\cal A} (d, \, -\alpha)
|x-y|^{-d-\alpha} \psi (m^{1/\alpha}|x-y|) =J(x,y)  \psi
(m^{1/\alpha}|x-y|)  $$ Thus the Dirichlet form ${\cal Q}$ of $X^m$
can also be written as follows
$$
{\cal Q}(u, v)\,=\, \int_{\R^d}\int_{\R^d} (u(x)-u(y))(v(x)-v(y))
J^m(x,y) dxdy.
$$

 For a
bounded $\kappa$-fat open set $D$, let  $X^{m,D}$ be killed
 relativistic stable process in $D$
with parameter $m > 0 $. From (\ref{DFC}), we have that the
Dirichlet form of $X^{m,D}$ is $({\cal Q}, {\cal F}^D)$, where $\cal
F ^D$ is given in (\ref{F^D}). For any $u, v\in {\cal F}^D$,
\begin{equation}\label{DF2}
{\cal Q}(u, v)=\int_D\int_D(u(x)-u(y))(v(x)-v(y))J^m(x, y)dxdy
+\int_Du(x)v(x)\kappa^m_D(x)dx,
\end{equation}
where
\begin{equation}\label{e:J}
\kappa^m_D(x):=2\int_{D^c}J^m(x, y)dy={\cal A} (d, \,
-\alpha)\int_{D^c}
 |x-y|^{-d-\alpha} \psi (m^{1/\alpha}|x-y|)dy .
\end{equation}

Let
$$
F^m(x,y):=\frac{J^m(x,y)}{J(x,y)}-1= \psi(m^{1/\alpha} |x-y|)-1$$
 and with $\kappa_D(x)$ given in (\ref{kappa_D})
$$
q^m(x):=\kappa^m_D(x)-\kappa_D(x) ={\cal A} (d, -\alpha) \int_{D^c}
F^m(x,y) |x-y|^{-d-\alpha} dy.
$$
Since $\psi(m^{1/\alpha} |x-y|)$ is a positive continuous function
and $D$ is bounded, $\inf_{x, y\in D} F^m(x, y) > -1$. We know
that $\psi(r)$ is a smooth function of $r^2$ and $\psi(0)=1$. Thus,
there exists a constant $c=c(D,\alpha, m)>0$ such that
\begin{equation}\label{Fbd}
 |F^m(x, y)| \leq c \,
|x-y|^2, \quad   \mbox{ for } x, y \in D.
\end{equation}
 Thus by Theorem \ref{t:A_2},
\begin{equation}\label{e:A}
F^m\in \A_2 (X^D) \quad\mbox{and} \quad \ln(1+F^m) \in \A_\infty
(X^D).
\end{equation}
Moreover, by (\ref{Fbd}),
$$
q^m(x) \le c_1 \int_{D^c} |x-y|^{-d-\alpha+2} dy \le c_2
\rho_D(x)^{2-\alpha}, \quad x \in D.
$$
So using the 3G theorem, $q^m \in S_\infty(X^D)$.

Let
\begin{eqnarray*}
&&K^m_t := \exp \left(\sum_{0<s \le t} \ln(1+F^m(X^D_{s-}, X^D_s)) -
{\cal A} (d, -\alpha)
\int^t_0 \int_D F^m(X^D_s, y)|X^D_s-y|^{-d-\alpha} dy ds\right.\\
&&\left.~~~~~~~~~~~~~~~~- \int_0^t q^m(X^D_s) ds\right).
\end{eqnarray*}
Using the multiplicative functional $K^m_t$, we define a semigroup
$Q^m_t$;
$$
Q^m_t f(x):=\E_x\left[f(X^D_t)K^m_t\right], \quad x \in D.
$$
Then by Theorem 4.8 of \cite{CS7}, $Q^m_t$ is a strongly continuous
semigroup in $L^2(D, dx)$ whose associated quadratic form is $({\cal
L}, {\cal F}^D)$, where
 \begin{eqnarray*}
{\cal L}(u, v)&=&\int_D\int_D(u(x)-u(y))(v(x)-v(y))(1+F^m(x,y))J(x, y)dxdy\\
&&+\int_Du(x)v(x)\kappa_D(x)dx + \int_Du(x)v(x) q^m(x)dx \\
&=&\int_D\int_D(u(x)-u(y))(v(x)-v(y))J^m(x, y)dxdy
+\int_Du(x)v(x)\kappa^m_D(x)dx
\end{eqnarray*}
for $u,v$ in ${\cal F}^D$. Thus by (\ref{DF2}),
 the quadratic form associated with $Q^m_t$
is exactly the Dirichlet form $({\cal Q}, {\cal F}^D)$ of $X^{m,D}$.
Therefore  $X^m$ can be obtained from $X$ through the Feynman-Kac
transform $K^m_t$. That is,
$$
\E_x\left[f(X^{m,D}_t)\right]= \E_x\left[f(X^D_t)K^m_t\right]
$$
for every positive Borel measurable function $f$. Since $X^{m,D}$ is
transient (for example, see Theorem 3.2 in \cite{CS4}), by Theorem
3.10 in \cite{C},
 we get the following
 theorem, which extends Theorem 3.1 in \cite{CS4}.
\medskip

\begin{thm}\label{T:6.0}
There exists Green function $V_D(x, y)$ for $X^{m}$ in $D$, jointly
continuous on $D\times D$ except along the diagonal, such that
$$ \int_D V_D(x, y) f(y) dy=\E_x \left[ \int_0^{\tau_D}
f(X^m_s) ds \right]
$$
for every Borel function $f\geq 0$  on $D$. Moreover, there exists
a positive constant $c=c(m,\alpha,D)$ such that
$$
c^{-1}\, G_D(x, y) \,\le \,V_D(x, y)\,\le\, c\,G_D(x, y), \quad x,y
\in D,
$$
and for every $x, y, z, w \in D$
\begin{equation} \label{e:VVV}
\frac{V_D(x,y) V_{D}(z,w)} { V_{D}(x,w)} \le c
\left(\frac{|x-w|\wedge |y-z|}{|x-y|  } \vee 1 \right)^\gamma
\left(\frac{|x-w|\wedge |y-z|}{|z-w|} \vee 1 \right)^\gamma
\frac{|x-w|^{d-\alpha}} {|x-y|^{d-\alpha} |z-w|^{d-\alpha}}.
\end{equation}
\end{thm}

\medskip

Since the jump measure of $X^m$ is $\frac12 v(x-y)$, we see that
L\'{e}vy system $(N^m, H^m)$ for $X^{m, D}$ can be chosen to be
\begin{equation}\label{e:Rm}
N^m(x,dy)=2J^m(x,y)dy=\frac{{\cal A}(d,-\alpha)\psi(m^{1/\alpha}
|x-y|)}{|x-y|^{d+\alpha}}dy \quad \mbox{in} \,\,D  \quad\quad
\hbox{and} \quad  H^m_t=t
\end{equation}
and the Revuz measure $\mu_{H^m} (dx)$ for $H^m$ is simply the
Lebesgue measure $dx$ on $D$ (cf. \cite{FOT}). Thus
\begin{equation}\label{e:Nm}
N^m(x,dy) \le N^D (x,dy)\quad \mbox{ in } D.
\end{equation}

\medskip

\begin{thm}\label{t:A_3}
If $D$ is a bounded $\kappa$-fat open set and $F$ is a function on
$D \times D$ with $|F(x,y)| \leq c|x-y|^{\beta}$ for some $\beta >
\alpha$ and $c > 0 $, then $F \in \A_2(X^{m,D})$ and
$$
\sup_{(x,w)\in D \times D} \E^w_x \left[ \sum_{s <
\tau_D}F(X^{m,D}_{s-}, X^{m,D}_s) \right] \, < \, \infty.
$$
\end{thm}
\pf By Theorem \ref{T:6.0} and (\ref{e:Nm}),
 $G_D(x,y) G_{D}(z,w) /G_{D}(x,w)$ is bounded above by
 (\ref{3G_est6}). So by (\ref{e:Nm}),
$$\frac {G_D(x,y)|F(y,z)|  G_D(z,w)}{G_D(x,w)}
\frac{N^m(y,dz)}{dz}
$$ is bounded above by
 (\ref{eqn:4.10}).
Now following the remainder argument in Theorem \ref{t:A_2}, we get
the result.
 \qed

\medskip

In Chen and Kim \cite{CK2}, an integral representation of
nonnegative excessive functions for the Schr\"odinger operator  is
established. Moreover it is shown that the Martin boundary is stable
under non-local Feynman-Kac perturbation. As mentioned in section 6
of \cite{CK2}, the method in \cite{CK2}
 works for a large class of strong Markov processes having a dual process.
 In particular,
Theorem 3.4 (3) \cite{CK2} is true for symmetric stable processes in bounded
$\kappa$-fat open sets. Applying Theorems 3.4 and 5.16 in \cite{CK2}
to symmetric stable and relativistic stable processes in bounded
$\kappa$-fat open sets, we have the following from our Theorem
\ref{T:6.0}.

\medskip

\begin{thm}\label{T:6.1}
For every $x\in D$ and $w\in \partial D$, $K(x, w):=
 \lim_{y\to w, y\in D} \frac{V_D(x, y)}{V_D(x_0, y)}$ exists and is finite.
It is jointly continuous on $D\times \partial D$. Moreover, there
exists  a positive constant $c$ such that
$$
c^{-1}\,M_D(x, w) \,\le \,K_D(x, w) \,\le \,c\,M_D(x, w), \quad x
\in D, w \in
\partial D,
$$ and for every $x, y, z, w \in D$
\begin{equation} \label{e:VKK}
\frac{V_D(x,y) K_{D}(z,w)} { K_{D}(x,w)} \le c
\left(\frac{|x-w|\wedge |y-z|}{|x-y|  } \vee 1 \right)^\gamma
\left(\frac{|x-w|\wedge |y-z|}{|z-w|} \vee 1 \right)^\gamma
\frac{|x-w|^{d-\alpha}} {|x-y|^{d-\alpha} |z-w|^{d-\alpha}}.
\end{equation}
\end{thm}

\begin{thm}\label{T:6.11}
 For every singular positive harmonic function $u$ for
$X^m$ in $D$, there is a unique finite measure $\nu$ on $\partial D$
such that
\begin{equation}\label{eqn:5.8}
u(x) =\int_{\partial D} K_D(x, z) \nu (dz).
\end{equation}
Thus the Martin boundary and the minimal Martin boundary for
relativistic $\alpha$-stable process $X^m$ in $D$ can all be
identified with the Euclidean boundary $\partial D$ of $D$.
\end{thm}

\medskip

When $u$ is singular harmonic in $D$ for $X^m$, the measure $\nu$ in
(\ref{eqn:5.8})  is called the Martin measure of $u$. Let
$H^+(X^{m,D})$ be the collection of all positive singular harmonic
function of $X^m$ in $D$. We denote $\E_x^h$ the expectation for the
conditional process obtained from  $X^{m,D}$ through Doob's
$h$-transform with $h(\cdot)$ starting from $x\in D$.

\medskip

\begin{prop}\label{c:conh2}
For every $\beta > \alpha$,
\begin{equation} \label{e:conh2}
\sup_{x\in D, h \in H^+(X^{m,D})} \E^h_x \left[ \sum_{s <
\tau_D}|X^{m,D}_{s-}-X^{m,D}_s|^{\beta} \right] \, < \, \infty.
\end{equation}
\end{prop}
\pf Using Theorem \ref{T:6.0}, \ref{T:6.1} and \ref{T:6.11}, the
proof is almost identical to the proof of Proposition \ref{c:conh}.
So we skip the proof here. \qed

\medskip

The boundary behavior of harmonic functions
 under non-local Feynman-Kac
perturbations has been studied by the first named author
\cite{K,K2}. In \cite{K2}, it is proven that  the ratio $u/h$ of two
singular harmonic functions for $X^m$ in an bounded $C^{1,1}$-open
set has non-tangential limits almost everywhere with respect to the
Martin measure of $h$. Due to (\ref{e:A}), now we can apply Theorem
4.7 in \cite{K2} to $X^m$ in  bounded $\kappa$-fat open sets and
extend the result in Theorem 4.11 \cite{K2}.

Recall  the Stolz open set for $\kappa$-fat open set $D$ from
\cite{K2}; For $Q \in \partial D$ and $\beta > (1-\kappa)/\kappa$,
let
$$
A^{\beta}_Q := \left\{ y \in D ;~ |y-Q| <  \beta \rho_D(y) \right\}.
$$

\begin{thm}\label{T:Fatou_R}
Suppose $D$ is a bounded $\kappa$-fat open set. If $k$ and $u$ are
singular harmonic functions for $X^m$ in $D$ and $\nu$ the Martin
measure for $k$. Then for $\nu$-a.e. $Q \in \partial D$,
$$
\lim_{ A^{\beta}_Q \ni x \rightarrow Q} \frac{u(x)}{k(x)}~ \mbox{
exists for every } \beta
>  \frac{1-\kappa}{\kappa} .
$$
\end{thm}

\section{Extension}

As we mentioned in the introduction of this paper, to keep our
exposition as transparent as possible, we have chosen rotationally
invariant symmetric $\alpha$-stable process to present our result.
However our approach works for a large class of  symmetric Markov
processes. We have not attempted to find the most general condition
where generalized 3G theorem  is true. However we like to point out
that the generalized 3G theorem extends to any symmetric Markov
processes satisfying the conditions below. We put these in terms of
Green function $G_D(x,y)$.

There exist positive constants $r_0$, $\gamma$, $M$ and $c_0$ such
that the following holds:
\begin{description}
\item{(C1)} (Lemma \ref{l2.5}) For all $Q\in \partial D$ and $r\in
(0, r_0)$, we have
$$
G_D(A_s(Q), z_0)\ge c_0(s/r)^{\gamma}G_D(A_r(Q), z_0), \qquad s\in
(0, r).
$$

\item{(C2)} (Theorem \ref{HP}) Let $y, x_{1}, x_{2}\in D$ and
$x_{1}, x_{2}\in D \setminus B(y, \rho_D(y)/2)$ such that
$|x_{1}-x_{2}|< L (\rho_D(x_1) \wedge \rho_D(x_2) )$. Then there
exists a constant $c:=c(L)$  such that $ G_D(x_{2},y)\,\leq\, c
G_D(x_{1}, y) .$

\item{(C3)} $G_D(x,y) \le c_0 |x-y|^{-d+\alpha}$ and  $G_D(x,y)
\ge c_0^{-1} |x-y|^{-d+\alpha}$ if  $|x-y| \le \rho_D(y)/2$.

\item{(C4)} (Lemma \ref{l:Green_L}) For every $Q \in
\partial D$ and  $r\in (0, r_0)$,
we have for $x, y \in D \setminus B(Q, r)$ and $z_1, z_2 \in D \cap
B(Q, r/M)$
$$
\frac{ G_{D}(x,z_1)}{ G_{D} (y,z_1)} \, \le \,c\, \frac{ G_{D}
(x,z_2) }{ G_{D} (y,z_2)}.
$$
\end{description}

It is well-known that (C2)-(C4) implies Green function estimates
(Theorem 3.4) (for example, see \cite{H} for a similar setup and the
conditions equivalent to the condition (C4)). Under the above
conditions (C1)-(C4), all the results in Sections 3 hold through the
same argument.

Moreover, if the symmetric Markov process is a stable (L\'{e}vy)
process, we know from the proof of Lemma \ref{l2.5} that
(\ref{e:low}) implies (C1) with $\gamma < \alpha$.

In the remainder of this section, we discuss the generalized 3G
theorem for another type of jump processes. Recently a class of (not
necessarily rotationally invariant) symmetric $\alpha$-stable
L\'{e}vy process was studied in \cite{BSS} and \cite{Sz2}. Let $
\sigma$ be the surface measure on $\partial B(0,1)={\{|y|= 1\}}$.
For any $\alpha\in (0, 2)$, $\alpha$-stable L\'{e}vy process
$Z=(Z_t, \P_x)$ is a symmetric L\'evy process in $\R^d$ such that
$$
\E\left[\exp(i\xi\cdot(Z_t-Z_0))\right]=\exp(-t \Phi(\xi)) \quad
\quad \mbox{ for every } x\in \R^d \mbox{ and } \xi\in \R^d.
$$
where
\begin{equation}\label{e:psi}
\Phi(\xi):= \int_{\{|y|= 1\}} |y \cdot \xi   |^{\alpha}
f(y)\sigma(dy)
\end{equation}
and $f$ is a symmetric function on ${\{|y|= 1\}}$ with $
0<c_1^{-1}\le f \le c_1< \infty$. When $f \equiv c$ for some
appropriate constant, $Z$ is a rotationally invariant symmetric
$\alpha$-stable process $X$, which we considered in Sections 2--3
(see \cite{BSS} and \cite{Sz2} for details).

Recall that a bounded domain $D$ is said to be Lipschitz if there is
a localization radius $r>0$  and a constant $\Lambda >0$ such that
for every $Q\in \partial D$, there is a Lipschitz  function $\phi_Q:
\R^{d-1}\to \R$ satisfying $\phi_Q (0)= 0$, $| \phi_Q (x)- \phi_Q
(z)| \leq \Lambda |x-z|$, and an orthonormal coordinate system
$y=(y_1, \cdots, y_{d-1}, y_d):=(\tilde y, y_d)$ such that $ B(Q,
r)\cap D=B(Q, r)\cap \{ y: y_d > \phi_Q (\tilde y) \}$. It is easy
to see that $D$ is $\kappa$-fat with characteristics $(r, \kappa)$
with $\kappa=(2 \sqrt{1+ \Lambda^2})^{-1}$.

>From \cite{Sz2}, we see that, if $D$ is  a bounded  Lipschitz, the
conditions (C2)--(C4) are true (see (2.3)--(2.4), Lemma 3.3--3.4 and
Theorem 4.2 in \cite{Sz2}). Moreover, by Lemma 3.3 in \cite{BSS},
the condition (C1) is true with $\gamma < \alpha$. Thus we have the
generalized 3G theorem  for $Z$.
\begin{thm}\label{33}
If $Z$ is symmetric $\alpha$-stable L\'{e}vy process and $D$ is
bounded  Lipschitz, the inequality (\ref{3G}) is true with some $0<
\gamma < \alpha$.
\end{thm}

Since the density of the L\'{e}vy measure for $Z$ is comparable to $
|x|^{-d-\alpha} $ ((2.1) in \cite{Sz2}), we have the following
theorem. The argument of the  proof is the same as the ones in the
proofs of Theorems \ref{t:A_2} and \ref{t:A_3}.  So we skip the
proof.

\begin{thm}\label{t:A_4} Suppose $Z$ is symmetric
$\alpha$-stable L\'{e}vy process and $D$ is  bounded  Lipschitz. If
$F$ is a function on $D \times D$ with $|F(x,y)| \leq
c|x-y|^{\beta}$ for some $\beta > \alpha$ and $c > 0 $, then $F \in
\A_2(Z^{D})$ and
$$
\sup_{(x,w)\in D \times D} \E^w_x \left[ \sum_{s <
\tau_D}F(Z^{D}_{s-}, Z^{D}_s) \right] \, < \, \infty.
$$
\end{thm}

{\it Acknowledgment}: The first named author of this paper thanks
Professor Renming Song and Professor Zhen-Qing Chen for their
encouragement and their helpful comments on the first version of
this paper.

\vspace{.5in}
\begin{singlespace}

\end{singlespace}
\end{doublespace}
\end{document}